\documentclass[10pt,a4paper]{article}
\usepackage[utf8]{inputenc}
\usepackage{amsthm,amsmath,enumerate,epic}
\usepackage[mathscr]{euscript}
\usepackage{parskip}
\usepackage{amssymb}
\usepackage{graphicx}
\usepackage{caption}
\usepackage{subcaption}
\usepackage{hyperref}
\usepackage[vcentermath]{youngtab}
\newtheorem{thm}{\bf Theorem}[subsection]
\newtheorem{prop}{Proposition}[subsection]

\newtheorem{remark}{\bf Remark}[thm]

\newtheorem{eg}{\bf Example}[subsection]
\newtheorem{theorem}{Theorem}[subsection]
\newtheorem{lemma}{\bf Lemma}[subsection]
\newtheorem{corollary}{\bf Corollary}[subsection]

\newcommand{\gtS}{\gtrdot}
\newcommand{\ltS}{\lessdot}
\newcommand{\bN}{\mathbb{N}}
\newcommand{\disjunion}{\,\,\dot{\cup}\,\,}
\newcommand{\la}[1]{\stackrel{#1}{\leftarrow}}

\newenvironment{Proof}{%
\par\noindent{\scshape Proof:}\begin{rm}}{\hfill$\Box$\end{rm}\newline}
\numberwithin{equation}{subsection}
\title{Local properties of Schubert varieties in the symplectic Grassmannian via a bounded RSK correspondence}
\date{}
\author {Shyamashree Upadhyay\\and\\Papi Ray\\Department of Mathematics\\ Indian Institute of Technology, Guwahati\\Assam-781039, INDIA\\emails: shyamashree@iitg.ac.in\\and\\popiroy93@iitg.ac.in}
\begin{document}
\maketitle
\begin{abstract}
In a paper by Ghorpade and Raghavan, they provide an explicit combinatorial description of the Hilbert function of the tangent cone at any point on a Schubert variety in the symplectic Grassmannian, by giving a certain ``degree-preserving'' bijection between a set of monomials defined by an initial ideal and a ``standard monomial basis". We prove here that this bijection is in fact a bounded RSK correspondence. 
\end{abstract}

\noindent \textbf{keywords:} Symplectic Grassmannian. Schubert variety. Tangent cone.  Hilbert function. RSK correspondence.

\section{Introduction}\label{s.Introduction}
The traditional RSK correspondence and its variants have many applications (see \cite{fulton}), particularly in combinatorics and representation theory. The main result of this paper is that the bijection of \cite{gr} is a variant of the RSK correspondence. This suggests a connection between the multiplicity problem for Schubert varieties and the applications of the RSK correspondence. It also provides a new perspective on the bijection of \cite{gr} and more evidence of the ubiquity of the RSK correspondence.\par
In \cite{kr} and \cite{kv-thesis}, an explicit Gr\"obner basis for the ideal of the tangent cone of a Schubert variety in the ordinary Grassmannian at a torus-fixed point was obtained. In \cite{kv}, the results of \cite{kr} and \cite{kv-thesis} were generalized to the case of Richardson varieties. In \cite{kv}, an explicit Gr\"obner basis for the ideal of the tangent cone at any $T$-fixed point of a Richardson variety in the ordinary Grassmannian was obtained, where $T$ denotes a maximal torus in the general linear group. In the study of singularities of Schubert varieties, Woo and Yong investigated Kazhdan-Lusztig ideals (see \cite{wy1}). These ideals encode coordinates and equations for neighborhoods of type $A$ Schubert varieties at torus fixed points. In \cite{wy2}, Woo and Yong provide a Gr\"obner basis for the Kazhdan-Lusztig ideals.\par
Sturmfels \cite{stur} and Herzog-Trung \cite{ht} proved results on a class of determinantal varieties which are equivalent to the results of \cite{kr},\cite{kv-thesis} and \cite{kv} for the case of Schubert varieties at the $T$-fixed point $e_{id}$. The key to their proofs was to use a version of the RSK correspondence (see \cite{fulton} for the classical RSK) in order to establish a degree-preserving bijection between a set of monomials defined by an initial ideal and a `standard monomial basis'. The proof given in \cite{kv} is based on a generalization of the RSK correspondence, which Kreiman calls the bounded RSK (BRSK). The BRSK is a degree-preserving bijection between an indexing set of a certain set of monomials defined by an initial ideal and an indexing set for a `standard monomial basis'. Result analogous to \cite{kv} has been obtained in \cite{su} for the orthogonal Grassmannian. In \cite{gk}, Graham and Kreiman give combinatorial descriptions of the restrictions to $T$-fixed points of the classes of structure sheaves of Schubert varieties in the $T$-equivariant $K$-theory of Grassmannians and of maximal isotropic Grassmannians of orthogonal and symplectic types. Graham and Kreiman also give formulas for the Hilbert series and Hilbert polynomials at $T$-fixed points of the corresponding Schubert varieties. These descriptions and formulas are given in terms of two equivalent combinatorial models: excited Young diagrams and set-valued tableaux.\par
In \cite{gr}, Ghorpade and Raghavan give an explicit combinatorial description of the multiplicity as well as the Hilbert function of the tangent cone at any point on a Schubert variety in the symplectic Grassmannian. In \cite{gr}, Ghorpade and Raghavan also provide an explicit Gr\"obner basis (with respect to certain term orders) for the ideal defining the tangent cone to a Schubert variety in the symplectic Grassmannian at any $T$-fixed point in it. The proof given in \cite{gr} also relies on a bijection between a set of monomials defined by an initial ideal and a `standard monomial basis' (see Chapter 6 of \cite{lr} for a standard monomial basis). In this paper, we prove that the bijection given in \cite{gr} is in fact a bounded RSK correspondence.\par
Schubert varieties are special cases of Richardson varieties. In \cite{kv}, Kreiman mentions in the introduction that the bijection given by him, when restricted to Schubert varieties in the ordinary Grassmannian is same as the bijection given in \cite{kr}. In this paper, our main theorem and its proof provides a proof of this fact.\par
The symplectic Grassmannian and a Schubert variety therein are defined in \S \ref{ss.notation} of this paper. In this paper, we consider the bijection given by Proposition 4.1 of \cite{gr} for Schubert varieties in the symplectic Grassmannian. As an application of our main theorem, we prove here that this bijection is a bounded RSK correspondence. The proof of this fact actually reduces to proving that the bijection $\tilde{\pi}$ given in \cite{kr} is the same as the map BRSK of \cite{kv}, the reduction being shown in \S \ref{ss.reduction} of this paper.\par
The organization of this paper is as follows. In \S \ref{s.maintheorem}, we state the main theorem (namely, Theorem \ref{t.main}), its corollary, and then we prove some important results needed to prove the main theorem. In \S \ref{s.theproof}, we provide a proof of the main theorem. The corollary to the main theorem (namely, Corollary \ref{c.main}) has been used in \cite{ru} to compute an explicit Gr\"obner basis for the ideal of the tangent cone at any $T$-fixed point of a Richardson variety in the symplectic Grassmannian, thus generalizing a result of Ghorpade and Raghavan.\par  
In \S \ref{s.theproblem}, we give an application of the main theorem after some necessary definitions and notation.
\section{The main theorem}\label{s.maintheorem}
In this section, we state the main theorem (Theorem \ref{t.main}) and its corollary. But before stating them, let us recall certain things from \cite{kr} and \cite{kv}, which are needed to state the main theorem.
\subsection{Recap}\label{ss.recap}
Given any positive integer $N$, we denote by $[N]$ the set $\{1,2,\ldots,N\}$. Given positive integers $d$ and $N$ with $d\leq N$, we denote by $I(d,N)$ the set of all $d$-element subsets of $[N]$. Let $\alpha=(\alpha_1, \ldots , \alpha_d)\in I(d,N)$, where $1\leq\alpha_1< \ldots < \alpha_d\leq N$. If $\beta=(\beta_1, \ldots ,\beta_d)\in I(d,N)$ be such that $1\leq\beta_1< \ldots <\beta_d\leq N$, then we say that $\alpha \leq \beta$ if $\alpha_i \leq \beta_i\ \forall i=1,\ldots,d$. Clearly, $\leq$ defines a partial order on $I(d,N)$.\par
Fix an element $v\in I(d,N)$. Let $\widetilde{\mathfrak{R}^v}$ denote the set of all ordered pairs $(r,c)$ such that $r\in[N]\setminus v$ and $c\in v$. Let $\widetilde{\mathfrak{N}^v}$ denote the subset of $\widetilde{\mathfrak{R}^v}$ consisting of those $(r,c)$ with $r>c$. Figure~\ref{figure.1} in Example~\ref{eg.etavtilde} shows the points of $\widetilde{\mathfrak{N}^v}$ for $v=(1,2,4,7,8,11)$. Let $\widetilde{T^v}$ denote the set of all monomials in $\widetilde{\mathfrak{N}^v}$. 
\begin{eg}\label{eg.etavtilde}
Let $d=6$, $N=13$, and $v=(1,2,4,7,8,11)$.
\vspace{.5 cm}

\begin{figure}[h]

\setlength{\unitlength}{0.3mm}   
\centering

\begin{picture}(60,100)(0,0)
\matrixput(0,0)(15,0){6}(0,15){7}{\circle*{1}}
\multiputlist(-15,0)(0,15){13,12,10,9,6,5,3}
\multiputlist(0,105)(15,0){1,2,4,7,8,11}

\linethickness{.2pt}
\dottedline{0.1}(0,90)(15,90)(15,75)(30,75)(30,45)(60,45)(60,15)(75,15)(75,0)

\put(75,5){$\longleftarrow$}
\put(85,5){ boundary of $\widetilde{\mathfrak{N}^v}$}
\end{picture}
\caption{The points of $\widetilde{\mathfrak{N}^v}$}
\label{figure.1}
\end{figure}
The points of the above grid represent the set $\widetilde{\mathfrak{R}^v}$ for $v=(1,2,4,7,8,11)$. The path sketched on the grid by some piece wise line segments denote the boundary of $\widetilde{\mathfrak{N}^v}$. The points on this grid which lie on the boundary of $\widetilde{\mathfrak{N}^v}$ or to the left of it are the points of $\widetilde{\mathfrak{N}^v}$.
\end{eg}
A \textbf{standard monomial} in $I(d,N)$ is a totally ordered sequence $\theta_1\geq\ldots\geq\theta_t$ of elements of $I(d,N)$. Such a standard monomial is called $v$-\textbf{compatible} if each $\theta_j$ is comparable to $v$ (with respect to the partial order $\leq$) but no $\theta_j$ equals $v$; it is \textbf{anti-dominated} by $v$ if $\theta_t\geq v$. Let $\widetilde{SM^{v,v}}$ denote the set of all $v$-compatible standard monomials in $I(d,N)$ anti-dominated by $v$.
\begin{eg}\label{e.smvvtilde}
Let $d=6$ and $N=13$. Then the totally ordered sequence $\theta_1=(4,7,8,11,12,13)\geq\theta_2=(3,5,6,9,10,11)\geq \theta_3=(1,2,3,6,7,10)$ is a standard monomial in $I(6,13)$. Let $v=(1,2,4,7,8,11)$. Then the above standard monomial (namely, $\theta_1\geq\theta_2\geq\theta_3$) is $v$-compatible because each $\theta_j\ (j=1,2,3)$ is comparable to $v$ but no $\theta_j$ equals $v$. In fact, in this case, we have $\theta_1\geq\theta_2\geq v$, $v\geq \theta_3$ and none of the $\theta_i$ ($i=1,2,3$) equals $v$.  Moreover, the standard monomial $\theta_1\geq\theta_2$ is anti-dominated by $v$ because $\theta_2\geq v$. The standard monomial $\theta_1\geq\theta_2$ is an element of $\widetilde{SM^{v,v}}$. 
\end{eg}
Given any ${\beta}_1=(r_1,c_1), \  {\beta}_2=(r_2,c_2)$ in $\widetilde{\mathfrak{N}^v}$, we say that ${\beta}_1=(r_1,c_1) > {\beta}_2=(r_2,c_2)$ if $r_1 >r_2$ and $c_1<c_2$. A sequence ${\beta}_1>\ldots>{\beta}_t$ of elements of $\widetilde{\mathfrak{N}^v}$ is called a $v$-\textbf{chain}. Given a $v$-chain ${\beta}_1=(r_1,c_1)>\ldots>{\beta}_t=(r_t,c_t)$, we define 
$$s_{\beta_1}\ldots s_{\beta_t}v:=(\{v_1,\ldots,v_d\}\setminus \{c_1,\ldots,c_t\})\cup \{r_1,\ldots ,r_t\}.$$
We say that an element $w$ of $I(d,N)$ \textbf{dominates} the $v$-chain ${\beta}_1>\ldots>{\beta}_t$ if $w \geq s_{\beta_1}\ldots s_{\beta_t}v $. Figure~\ref{figure.2} in Example~\ref{eg.vchain} shows a $v$-chain for $v=(1,2,4,7,8,11)$.
\begin{eg}\label{eg.vchain}
Let $d=6$, $N=13$, and $v=(1,2,4,7,8,11)$.\\
\vspace{.5 cm}
\begin{figure}[h]
\setlength{\unitlength}{0.3mm}
\centering
\begin{picture}(60,100)(0,0)
\matrixput(0,0)(15,0){6}(0,15){7}{\circle*{.9}}
\multiputlist(-15,0)(0,15){13,12,10,9,6,5,3}
\multiputlist(0,105)(15,0){1,2,4,7,8,11}
\put(0,0){\makebox(0,0){$\bullet$}}
\put(15,15){\makebox(0,0){$\bullet$}}
\put(30,30){\makebox(0,0){$\bullet$}}
\put(60,45){\makebox(0,0){$\bullet$}}
\linethickness{.2pt}
\dottedline{0.1}(0,90)(15,90)(15,75)(30,75)(30,45)(60,45)(60,15)(75,15)(75,0)
\end{picture}
\caption{A $v$-chain}
\label{figure.2}
\end{figure}
The four dark circles in this grid denote a $v$-chain given by $\beta_1 > \beta_2 > \beta_3 > \beta_4$, where $\beta_1 =(13,1)$, $\beta_2 =(12,2)$, $\beta_3  =(10,4)$, $\beta_4=(9,8)$. If we take $w=(8,9,10,11,12,13)$, then this $v$-chain is dominated by $w$ because $s_{\beta_1}\ldots s_{\beta_4}v=(7,9,10,11,12,13)\leq w$.
\end{eg}
Let $\mathfrak{S}$ be a monomial in $\widetilde{\mathfrak{R}^v}$. By a $v$-chain in $\mathfrak{S}$, we mean a sequence ${\beta_1}>\ldots>{\beta_t}$ of elements of $\mathfrak{S}\cap \widetilde{\mathfrak{N}^v}$. We say that $w$ \textbf{dominates} $\mathfrak{S}$ if $w$ dominates every $v$-chain in $\mathfrak{S}$.\par 
We call \textbf{distinguished} the subsets $\mathfrak{S}$ of $\widetilde{\mathfrak{N}^v}$ satisfying the following conditions:

\noindent (A) For $(r,c)\neq (r',c')$ in $\mathfrak{S}$, we have $r\neq r'$ and $c\neq c'$.\\
\noindent (B) If $\mathfrak{S}=\{(r_1,c_1),\ldots,(r_p,c_p)\}$ with $r_1<r_2<\ldots<r_p$, then for $j$, $1\leq j\leq p-1$, we have either $c_j>c_{j+1}$ or $r_j<c_{j+1}$.
\begin{eg}\label{e.distinguished}
For $v=(1,2,4,7,8,11)$, the subset $\mathfrak{S}$ of $\widetilde{\mathfrak{N}^v}$ given by $\mathfrak{S}=\{(3,2),(5,1),(9,8),(12,7)\}$ is distinguished.
\end{eg}
\begin{remark}\label{r.distinguished}
By Proposition 4.3 of \cite{kr}, there exists a bijection between elements $w$ of $I(d,N)$ satisfying $w\geq v$ on the one hand and distinguished subsets of $\widetilde{\mathfrak{N}^v}$ on the other hand. We denote this bijective correspondence by $w\leftrightarrow\mathfrak{S}_w$. 
\end{remark}
Let $\mathfrak{S}$ be a non-empty monomial in $\widetilde{\mathfrak{N}^v}$. If $\beta_1>\cdots>\beta_t$ is a $v$-chain in $\mathfrak{S}$, then we call $\beta_1$ the \textbf{head} of the $v$-chain and $\beta_t$ its \textbf{tail}. We call $t$ to be the \textbf{length} of the $v$-chain. We say that an element $\beta$ of $\mathfrak{S}$ is $t$-\textbf{deep} in $\mathfrak{S}$ (where $t$ is a positive integer) if $\beta$ is the tail of a $v$-chain in $\mathfrak{S}$ of length $t$. The \textbf{depth} of $\beta$ in $\mathfrak{S}$ is defined to be $t$ if $\beta$ is $t$-deep in $\mathfrak{S}$ but not $(t+1)$-deep in $\mathfrak{S}$.
\begin{eg}\label{e.depth}
Let $v=(1,2,4,7,8,11)$ and $\mathfrak{S}=\{(9,1),(6,2),(5,4),(13,8),(12,11)\}$ be a monomial in $\widetilde{\mathfrak{N}^v}$. Let $\beta=(5,4)$. Then it is easy to see that $\beta$ is $1$-deep, $2$-deep, and $3$-deep in $\mathfrak{S}$. But $\beta$ is not $4$-deep in $\mathfrak{S}$. In fact, $(9,1)>(6,2)>(5,4)$ is a $v$-chain in $\mathfrak{S}$ and this is a $v$-chain in $\mathfrak{S}$ of maximum length having $\beta=(5,4)$ as its tail. Hence the depth of $\beta$ in $\mathfrak{S}$ is $3$ here.
\end{eg}
We will now recall the map $\pi$ of \cite{kr}. Let $\mathfrak{S}$ be a non-empty monomial in the elements of $\widetilde{\mathfrak{N}^v}$. We partition $\mathfrak{S}$ in two stages. First we partition $\mathfrak{S}$ into subsets $\mathfrak{S}_1,\ldots,\mathfrak{S}_k$, where $k$ is the largest length of a $v$-chain in $\mathfrak{S}$: $\beta\in\mathfrak{S}$ belongs to $\mathfrak{S}_j$ if it is $j$-deep but not $(j+1)$-deep.\par
Now we partition each $\mathfrak{S}_j$ into subsets called \textit{blocks} as follows. We arrange the elements of $\mathfrak{S}_j$ in non-decreasing order of their row numbers (all arrangements are from left to right; and elements occur with their respective multiplicities). Among those with the same row number, the arrangement is by non-decreasing order of column numbers. Two consecutive members $(r,c)$, $(R,C)$ in this arrangement are said to be \textit{related} if $r>C$. The blocks are the equivalence classes of the smallest equivalence relation containing the above relations.\par
Let $\mathfrak{B}$ be a single block of some $\mathfrak{S}_j$. Let
$$(r_1,c_1),\ldots,(r_p,c_p)$$
be the elements of $\mathfrak{B}$ written in non-decreasing order of both row and column numbers (in such an arrangement, the elements occur with their respective multiplicities). We set $w(\mathfrak{B}):=(r_p,c_1)$ and $\mathfrak{B}'$ to be the monomial
$$\{(r_1,c_2),(r_2,c_3),\ldots,(r_{p-2},c_{p-1}),(r_{p-1},c_p)\}.$$
Set $\mathfrak{S}_j^{(1)}:=\cup_{\mathfrak{B}}\mathfrak{B}'$ (where the index $\mathfrak{B}$ runs over all blocks of $\mathfrak{S}_j$) and $\mathfrak{S}^{(1)}:=\cup_{j=1}^{k}\mathfrak{S}_j^{(1)}$. It follows from Corollary 4.13 of \cite{kr} that the set
$$\{w(\mathfrak{B})|\mathfrak{B}\ \text{is a block of}\ \mathfrak{S}\}$$
is a distinguished subset of $\widetilde{\mathfrak{N}^v}$. Let $w$ be the corresponding element of $I(d,N)$ (under the correspondence given in Remark \ref{r.distinguished}). Set 
$$\pi(\mathfrak{S}):=(w,\mathfrak{S}^{(1)}).$$
This finishes the description of the map $\pi$ of \cite{kr}. Example~\ref{eg.themappi} below gives a detailed illustration of the map $\pi$ of \cite{kr}.
\begin{eg}\label{eg.themappi}
Let $d=6$, $N=13$, and $v=(1,2,4,7,8,11)$. The dark circles in the grid in Figure~\ref{figure.3} represent a monomial $\mathfrak{S}$ in $\widetilde{\mathfrak{N}^v}$, where $\mathfrak{S}=\{(3,2),(5,4),(6,2),(9,1),(9,1),(10,7),(10,7),(10,7),(10,8),(12,1),(13,4)\}$. The numbers written near the dark circles denote the multiplicities of these elements in the monomial $\mathfrak{S}$. For this monomial $\mathfrak{S}$, we have
$$\mathfrak{S}_1=\{(9,1),(9,1),(12,1),(13,4)\}$$
$$\mathfrak{S}_2=\{(3,2),(6,2),(10,7),(10,7),(10,7),(10,8)\}$$
$$\mathfrak{S}_3=\{(5,4)\}$$
Here $\mathfrak{S}_1$ and $\mathfrak{S}_3$ are single blocks. And $\mathfrak{S}_2$ has two blocks given by $\{(3,2),(6,2)\}$ and $\{(10,7),(10,7),(10,7),(10,8)\}$. The dark line segments on the grid show the block decomposition of the monomial $\mathfrak{S}$. The set 
$$\{w(\mathfrak{B})|\mathfrak{B}\ \mbox{is a block of}\ \mathfrak{S}\}=\{(13,1),(6,2),(10,7),(5,4)\}.$$
\vspace{.5 cm}
\begin{figure}[h]
\setlength{\unitlength}{0.6mm}
 \centering
\begin{picture}(60,100)(0,0)
\matrixput(0,0)(15,0){6}(0,15){7}{\circle*{.9}}
\multiputlist(-15,0)(0,15){13,12,10,9,6,5,3}
\multiputlist(0,105)(15,0){1,2,4,7,8,11}
\put(0,15){\makebox(0,0){$\bullet$}}
\put(0,45){\makebox(0,0){$\bullet$}}
\put(15,60){\makebox(0,0){$\bullet$}}
\put(15,90){\makebox(0,0){$\bullet$}}
\put(30,0){\makebox(0,0){$\bullet$}}
\put(30,75){\makebox(0,0){$\bullet$}}
\put(45,30){\makebox(0,0){$\bullet$}}
\put(60,30){\makebox(0,0){$\bullet$}}
\linethickness{.7pt}
\dottedline{0.1}(0,45)(0,0)(30,0)
\dottedline{0.1}(45,30)(60,30)
\dottedline{0.1}(15,90)(15,60)
\put(2,17){1}
\put(2,47){2}
\put(17,62){1}
\put(17,92){1}
\put(32,2){1}
\put(32,77){1}
\put(47,32){3}
\put(62,32){1}
\linethickness{.2pt}
\dottedline{0.1}(0,90)(15,90)(15,75)(30,75)(30,45)(60,45)(60,15)(75,15)(75,0)
\put(60,20){$\longleftarrow$}
\put(70,20){ boundary of $\widetilde{\mathfrak{N}^v}$}
\end{picture}
\caption{A monomial $\mathfrak{S}$}
\label{figure.3}
\end{figure}
Therefore $w=(5,6,8,10,11,13)$, and  $\mathfrak{S}^{(1)}=\{(9,1),(9,1),(12,4),(3,2),(10,7),(10,7),(10,8)\}$.
\end{eg}
Using $\pi$, we now recall the  map $\tilde{\pi}$ of \cite{kr} from $\widetilde{T^v}$ to $\widetilde{SM^{v,v}}$. Proceed by induction on the degree of an element $\mathfrak{S}$ of $\widetilde{T^v}$. The image of the empty monomial under $\tilde{\pi}$ is taken to be the empty monomial. Let $\mathfrak{S}$ be non-empty, and suppose that $\pi(\mathfrak{S}) =(w,\mathfrak{S}^{(1)})$. By (1) and
(2) of Proposition~4.1 of \cite{kr}, the degree of $\mathfrak{S}^{(1)}$ is strictly less than that of $\mathfrak{S}$,  and so by induction $\tilde{\pi}(\mathfrak{S}^{(1)})$ is defined.
Suppose that $\tilde{\pi}(\mathfrak{S}^{(1)})=w'\geq\ldots$. By induction we also know that the degree of $\mathfrak{S}^{(1)}$ is the same as that of $w'\geq\ldots$ and that $w'$ is the least element of $I(d,N)$ that dominates $\mathfrak{S}^{(1)}$. By (3) of
Proposition~4.1 of \cite{kr}, we have $w\geq w'$, and we set
$\tilde{\pi}(\mathfrak{S}):=w\geq\tilde{\pi}(\mathfrak{S}^{(1)})$. This finishes the description of the map $\tilde{\pi}$ of \cite{kr}.
\begin{eg}\label{eg.tildepi}
For the monomial $\mathfrak{S}$ in Example~\ref{eg.themappi} above, we have $\tilde{\pi}(\mathfrak{S})=(5,6,8,10,11,13)\geq (3,4,8,10,11,12)\geq (2,4,7,8,10,11)\geq (1,2,7,8,10,11)\geq (1,2,4,8,9,11)\geq (1,2,4,7,9,11)$.
\end{eg}
Let us now recall the map $BRSK$ of \cite{kv}. But for recalling the map $BRSK$ of \cite{kv}, we first need to recall a lot of definitions from \cite{kv}.\par
A \textbf{Young diagram} (resp. \textbf{notched diagram}) is a collection of boxes arranged into a left and top justified array (resp. into left justified rows). The \textbf{empty Young diagram} is the Young diagram with no boxes. A notched diagram may contain rows with no boxes; however, a Young diagram may not, unless it is the empty Young diagram. A \textbf{Young tableau} (resp. \textbf{notched tableau}) is a filling of the boxes of a Young diagram (resp. notched diagram) with positive integers. The \textbf{empty Young tableau} is the Young tableau with no boxes. Let $P$ be either a notched tableau or a Young tableau. We say that $P$ is \textbf{row strict} if the entries of any row of $P$ strictly increase as one moves to the right. If $P$ is a Young tableau, then we say that $P$ is \textbf{semistandard} if it is row strict and the entries of any column weakly increase as one moves down.
\begin{eg}\label{eg.rowstrict-ss}
A row strict notched tableau $P$, and a semistandard Young tableau $R$.
$$P=\young(1234,24,567,4578),\ \ R=\young(1234,2456,35,6)$$
\end{eg}
Let $P$ be a row strict notched tableau, and $b$ be a positive integer. Since $P$ is row strict, its entries which are greater than or equal to $b$ are right justified in each row. If we remove these entries (which are greater than or equal to $b$) from $P$, then we are left with a row strict notched tableau, which we denote by $P^{<b}$. We say that $P$ is \textbf{semistandard on} $b$ if $P^{<b}$ is a semistandard Young tableau. 
\begin{eg}\label{eg.P<b}
For the row strict notched tableau $P$ in Example~\ref{eg.rowstrict-ss}, and $b=4$, we have
$$P^{<b}=\young(123,2).$$
However, for the same $P$, if we take $b=6$, then
$$P^{<b}=\young(1234,24,5,45).$$
Hence $P$ is semistandard on $4$, but not on $6$.
\end{eg}
Let us now recall the \textbf{ordinary Schensted insertion} process from \S 3 of \cite{kv}. It is an algorithm which takes as input a semistandard Young tableau $P$ and a positive integer $a$, and produces as output a new semistandard Young tableau with the same shape as $P$ plus one extra box, and with the same entries as $P$ (possibly in different locations) plus one additional entry, namely $a$. To begin, insert $a$ into the first row of $P$, as follows. If $a$ is strictly bigger than all entries in the first row of $P$, then place $a$ in a new box on the right end of the first row, and the insertion process terminates. Otherwise, find the smallest entry of the first row of $P$ which is greater than or equal to $a$, and replace that number with $a$. We say that the number which was replaced was ``bumped" from the first row. Insert the bumped number into the second row in precisely the same way as $a$ was inserted into the first row. This process continues down the rows until, at some point, a number is placed in a new box on the right end of some row, at which point the process terminates.\par
We next describe the \textbf{bounded insertion algorithm}, which
takes as input a positive integer $b$, a notched tableau $P$ which
is semistandard on $b$, and a positive integer $a<b$, and produces
as output a notched tableau which is semistandard on $b$, which we
denote by $P\la{b}a$.
\begin{quote}
\noindent\textbf{Bounded Insertion} {\it
\begin{itemize}
\item[\textbf{Step 1.}] Remove all entries of $P$ which are
greater than or equal to $b$ from $P$, resulting in the
semistandard Young tableau $P^{<b}$.
\item[\textbf{Step 2.}] Insert $a$ into $P^{<b}$ using the
ordinary Schensted insertion process (as described above).
\item[\textbf{Step 3.}] Place the entries of $P$ which were
removed when forming $P^{<b}$ in Step 1 back into the Young
tableau resulting from Step 2, in the same rows from which they
were removed.
\end{itemize}
}
\end{quote}
\begin{eg}\label{eg.bounded-insertion}
Let $P=\young(1234,24,567,4578)$, $a=3$ and $b=4$. We compute $P\la{b}a$. Observe that in Step 1, we obtain $P^{<4}=\young(123,2)$. In Step 2, we insert $a=3$ into $P^{<4}$ using the ordinary Schensted insertion process, to get $\young(123,23)$. And finally in Step 3, we obtain $P\la{4}3=\young(1234,234,567,4578)$.
\end{eg}
Let $S$ be any set.  A \textbf{multiset} $E$ on $S$ is defined to be a function $E:S\rightarrow \{0,1,2,\ldots,\}$.  One should think of $E$ as consisting of the set $S$ of elements, but with each $s\in S$ occurring $E(s)$ times. Note that a set is a special type of multiset in which each element occurs exactly once. We call $E(s)$ the \textbf{degree} or \textbf{multiplicity} of $s$
in $E$. Define the multiset $E\disjunion F$ as follows:
\begin{align*}
(E\disjunion F)(s)&=E(s)+F(s),\ s\in S
\end{align*}
Let $\bN$ denote the positive integers. Let $A=\{a_1,a_2,\ldots\}$ and $B=\{b_1,b_2,\ldots\}$ be two multisets on $\bN$ of the same degree, with $a_i\leq a_{i+1}$, $b_i\leq b_{i+1}$, for all $i$. We say that $A$ is less than or equal to $B$ in the \textbf{termwise order} if $a_i\leq b_i$ for all $i$. We denote this by $A\leq B$. We say that $A$ is less than $B$ in the \textbf{strict termwise order} if $a_i< b_i$ for all $i$.  We denote this by $A\ltS B$.\par If $A$, $B$, $C$, and $D$ are multisets on $\bN$ such that
$|A\disjunion D|=|B\disjunion C|$, then we write
\begin{equation}\label{e.m.set_subtraction}
A-C\leq B-D \hbox{ to indicate that }A\disjunion D\leq B\disjunion
C.
\end{equation}
Let  $U=\{(e_1,f_1),(e_2,f_2),\ldots\}$ be a multiset on $\bN^2$.
Define $U_{(1)}$ and $U_{(2)}$ to be the multisets
$\{e_1,e_2,\ldots\}$ and $\{f_1,f_2,\ldots\}$ respectively on
$\bN$. Define the \textbf{nonvanishing}, \textbf{negative}, and
\textbf{positive parts} of $U$ to be the following multisets:
\begin{align*}
U^{\neq 0}&=\{(e_i,f_i)\in U\mid e_i-f_i\neq 0\},\\
U^-&=\{(e_i,f_i)\in U\mid e_i-f_i<0\},\\
U^+&=\{(e_i,f_i)\in U\mid e_i-f_i>0\}.
\end{align*}
We say that $U$ is \textbf{nonvanishing} if $U\subset (\mathbb{N}^2)^{\neq
0}$, \textbf{negative} if $U\subset (\mathbb{N}^2)^-$, and
\textbf{positive} if $U\subset (\mathbb{N}^2)^+$.  Impose the following transitive relation on multisets on $\bN^2$:
\begin{equation}
U\leq V\iff U_{(1)}-U_{(2)}\leq V_{(1)}-V_{(2)}.
\end{equation}
Let $\iota$ be the map on multisets on $\mathbb{N}^2$ defined by, 
$$\iota(\{(e_1,f_1),(e_2,f_2),\ldots\}):=\{(f_1,e_1),(f_2,e_2),\ldots\}.$$
Then $\iota$ is an involution, and it maps negative multisets on $\mathbb{N}^2$ to positive ones and visa-versa.\par
A \textbf{notched bitableau} is a pair $(P, Q)$ of notched
tableaux of the same shape (i.e., the same number of rows and the
same number of boxes in each row).  The \textbf{degree} of $(P,Q)$
is the number of boxes in $P$ (or $Q$). A notched bitableau
$(P,Q)$ is said to be \textbf{row strict} if both $P$ and $Q$ are
row strict. A row strict notched bitableau $(P,Q)$ is said to be
\textbf{semistandard} if
\begin{equation}\label{e.s.ss_tabl_1}
P_1-Q_1\leq\cdots\leq P_r-Q_r,
\end{equation}
where $r$ is the total number of rows in $P$ (or $Q$) and for each $i\in\{1,\ldots,r\}$, $P_i$ (resp. $Q_i$) denotes the $i$-th row (from the top) of $P$ (resp. $Q$). A row strict notched bitableau $(P,Q)$ is said to be
\textbf{negative} if $P_i\ltS Q_i$, $i=1,\ldots,r$,
\textbf{positive} if $P_i\gtS Q_i$, $i=1,\ldots,r$, and
\textbf{nonvanishing} if either
\begin{equation}\label{e.s.ss_tabl_2}
P_i\ltS Q_i\ \ \text{   or   }\ \ P_i\gtS Q_i,
\end{equation}
for each $i=1,\ldots,r$.
\begin{eg}\label{eg.ss-notched-bitableau}
Consider the notched bitableau
$$(P,Q)=\left(\ \young(123,4567)\ \ \ ,
\ \ \ \young(789,2345)\right).$$ We have that
\begin{itemize}
\item[1.] $(P,Q)$ is row strict.

\item[2.] $P_1\disjunion Q_2=\{1,2,2,3,3,4,5\}\leq
\{4,5,6,7,7,8,9\}=P_2\disjunion Q_1$. Therefore, $P_1-Q_1\leq
P_2-Q_2$. Thus $(P,Q)$ is semistandard.

\item[3.] $P_1\ltS Q_1$, and $P_2\gtS Q_2$. Thus $(P,Q)$ is nonvanishing.
\end{itemize}
\end{eg}
We next define the \textbf{bounded RSK correspondence}, $BRSK$, a
function which maps negative multisets on $\mathbb{N}^2$ to negative
semistandard notched bitableaux. Let
$U=\{(a_1,b_1),\ldots,(a_t,b_t)\}$ be a negative multiset on
$\mathbb{N}^2$, whose entries we assume are listed in
\textbf{lexicographic order}: (i) $b_1\geq\cdots\geq b_t$, and
(ii) if for any $i\in\{1,\ldots,t-1\}$, $b_i=b_{i+1}$, then
$a_i\geq a_{i+1}$. We inductively form a sequence of notched
bitableaux $(P^{(0)},Q^{(0)})$, $(P^{(1)},Q^{(1)}),$
$\ldots,(P^{(t)},Q^{(t)})$, such that $P^{(i)}$ is semistandard on
$b_i$, $i=1,\ldots,t$, as follows:
\begin{quote}
Let $(P^{(0)},Q^{(0)})=(\emptyset,\emptyset)$, and let $b_0=b_1$. Assume
inductively that we have formed $(P^{(i)}, Q^{(i)})$, such that
$P^{(i)}$ is semistandard on $b_i$, and thus on $b_{i+1}$, since
$b_{i+1}\leq b_i$. Define $P^{(i+1)}=P^{(i)}\la{b_{i+1}}a_{i+1}$.
Since bounded insertion preserves semistandardness on $b_{i+1}$,
$P^{(i+1)}$ is also semistandard on $b_{i+1}$. Let $j$ be the row
number of the new box of this bounded insertion. Define
$Q^{(i+1)}$ to be the notched tableau obtained by placing
$b_{i+1}$ on the {\it left} end of row $j$ of $Q^{(i)}$ (and
shifting all other entries of $Q^{(i)}$ to the {\it right} one box).
Clearly $P^{(i+1)}$ and $Q^{(i+1)}$ have the same shape.
\end{quote}
Then $BRSK(U)$ is defined to be $(P^{(t)},Q^{(t)})$. In the process above, we write $(P^{(i+1)},Q^{(i+1)})=(P^{(i)},Q^{(i)})\la{b_{i+1}}a_{i+1}$. In terms of this notation,
$$BRSK(U)=((\emptyset,\emptyset)\la{b_1}a_1)\cdots\la{b_t}a_t.$$
\begin{eg}\label{eg.brsk}
Let $U=\{(3,4),(4,6),(4,7),(5,7),(2,4)\}$ be a negative multiset on $\mathbb{N}^2$. Arranging $U$ in lexicographic order, we have $U=\{(5,7),(4,7),(4,6),(3,4),(2,4)\}$. Then
\[
\begin{array}{l@{\hspace{.9cm}}l}
P^{(0)}=\emptyset 
&Q^{(0)}=\emptyset\\ \\
P^{(1)}=\emptyset\la{7}5=\young(5)
&Q^{(1)}=\young(7)\\ \\
P^{(2)}=\young(5)\la{7}4=\young(4,5)
&Q^{(2)}=\young(7,7)\\ \\
P^{(3)}=\young(4,5)\la{6}4=\young(4,4,5)
&Q^{(3)}=\young(7,7,6)\\ \\
P^{(4)}=\young(4,4,5)\la{4}3=\young(34,4,5)
&Q^{(4)}=\young(47,7,6)\\ \\
P^{(5)}=\young(34,4,5)\la{4}2=\young(24,34,5)
&Q^{(5)}=\young(47,47,6)\\
\end{array}
\]
\vspace{.5em}
Therefore $BRSK(U)=\left(\, \young(24,34,5)\ ,\
\young(47,47,6)\right)$.
\end{eg}
If $(P,Q)$ is a non-vanishing semistandard notched bitableau, we define $\iota(P,Q)$ to be the notched bitableau obtained by reversing the order of the rows of $(Q,P)$. If $U$ is a positive multiset on $\mathbb{N}^2$, then $BRSK(U)$ is defined to be $\iota(BRSK(\iota(U)))$.
\subsection{Statement of the main theorem}\label{ss.tmain}
The main results of this paper are Theorem~\ref{t.main} and Corollary~\ref{c.main} below.
\begin{theorem}\label{t.main}
Let $U$ be a finite monomial in $\widetilde{\mathfrak{N}^v}$. Let $\pi (U)=(w_0,U^{(1)})$, $\pi(U^{(1)})=(w_1,U^{(2)}),\ldots$ and so on till $\pi(U^{(m)})=(w_m,\emptyset)$, where $\emptyset$ is the empty monomial . Then for each $r \in \{0,1,\ldots,m\}$, the following holds :\\
(i) All the row numbers of the distinguished subset $\mathfrak{S}_{w_r}$ (corresponding to $w_r$) consist of the $((m+1)-r)$-th row entries of the left hand notched tableau of $BRSK(U)$.\\
(ii) All the column numbers of $\mathfrak{S}_{w_r}$ comprise of the $((m+1)-r)$-th row entries of the right hand notched tableau of $BRSK(U)$.
\end{theorem}
Example~\ref{eg.maintheorem} below illustrates Theorem~\ref{t.main}.
\begin{eg}\label{eg.maintheorem}
Let $d=6$, $N=13$, and $v=(1,2,4,7,8,11)$. Let $U$ be a finite monomial in $\widetilde{\mathfrak{N}^v}$ given by
$$U=\{(3,2),(5,4),(6,2),(9,1),(9,1),(10,7),(10,7),(10,7),(10,8),(12,1),(13,4)\}.$$
Then $\pi(U)=(w_0,U^{(1)})$, $\pi(U^{(1)})=(w_1,U^{(2)})$, $\ldots$, and so on till $\pi(U^{(5)})=(w_5,\emptyset)$, where
$$w_0=(5,6,8,10,11,13),w_1=(3,4,8,10,11,12),w_2=(2,4,7,8,10,11),$$
$$w_3=(1,2,7,8,10,11),w_4=(1,2,4,8,9,11),w_5=(1,2,4,7,9,11)$$
and
$$U^{(1)}=\{(9,1),(9,1),(12,4),(3,2),(10,7),(10,7),(10,8)\},$$
$$U^{(2)}=\{(9,1),(9,4),(10,7),(10,8)\},U^{(3)}=\{(9,4),(9,7),(10,8)\},$$
$$U^{(4)}=\{(9,7),(9,8)\},U^{(5)}=\{(9,8)\}.$$
\newcommand{\J}{10}
\newcommand{\B}{12}
\newcommand{\A}{13}
Here $m=5$ and $BRSK(U)=\left(\hspace{0.2cm}\young(9,9,\J,\J,3\J\B,56\J\A),\hspace{0.5cm}\young(8,7,4,1,127,1247)\hspace{0.2cm}\right)$.\\ \\
Also, $\mathfrak{S}_{w_0}=\{(13,1),(6,2),(10,7),(5,4)\}$, $\mathfrak{S}_{w_1}=\{(12,1),(3,2),(10,7)\}$,  $\mathfrak{S}_{w_2}=\{(10,1)\}$,   $\mathfrak{S}_{w_3}=\{(10,4)\}$,  $\mathfrak{S}_{w_4}=\{(9,7)\}$ and $\mathfrak{S}_{w_5}=\{(9,8)\}$. For each $r\in\{0,1,\ldots,5\}$, conditions (i) and (ii) of the above theorem hold true. For example, if we take $r=1$, then all the row numbers of the distinguished subset $\mathfrak{S}_{w_1}$ are given by $\{3,10,12\}$, which are the $((5+1)-1)$-th row entries of the left hand notched tableau of $BRSK(U)$. Similarly, for the column numbers of $\mathfrak{S}_{w_1}$, which comprise of the $((5+1)-1)$-th row entries of the right hand notched tableau of $BRSK(U)$.  
\end{eg}
\begin{corollary}\label{c.main}
For any finite monomial $U$ in $\widetilde{\mathfrak{N}^v}$, $\tilde{\pi}(U)=BRSK(U)$. 
\end{corollary}
\begin{remark}\label{r.c.main}
In Corollary \ref{c.main} above, when we say that $\tilde{\pi}(U)=BRSK(U)$, we mean the following: Let $BRSK(U)=(P,Q)$. Let $P_1,\ldots,P_r$ (resp. $Q_1,\ldots,Q_r$) denote all the rows of $P$ (resp. $Q$) from top to bottom. Then we mean that $\tilde{\pi}(U)$ equals the standard monomial $P_r-Q_r\geq\cdots\geq P_1-Q_1$ (where $P_i-Q_i$ has to be interpreted as $P_i\disjunion(v\setminus Q_i)$ in the sense of \S 4 of \cite{kv}).
\end{remark}
Example~\ref{eg.maincorollary} below illustrates Corollary~\ref{c.main}.
\begin{eg}\label{eg.maincorollary}
For the monomial $U$ in Example~\ref{eg.maintheorem}, we have $\tilde{\pi}(U)=w_0\geq w_1\geq\cdots\geq w_5$, where $w_0,\ldots,w_5$ are as given in Example~\ref{eg.maintheorem}. It is now easy to verify (using Remark~\ref{r.c.main} above) that $\tilde{\pi}(U)=BRSK(U)$. $BRSK(U)$ is given in Example~\ref{eg.maintheorem} above. Let $P_1,\ldots,P_6$ denote the rows (from top to bottom) of the left hand notched tableau of $BRSK(U)$. Similarly, let $Q_1,\ldots,Q_6$ denote the rows (from top to bottom) of the right hand notched tableau of $BRSK(U)$. Then observe that $w_0=P_6-Q_6,w_1=P_5-Q_5,w_2=P_4-Q_4,w_3=P_3-Q_3,w_4=P_2-Q_2$ and $w_5=P_1-Q_1$. Hence $\tilde{\pi}(U)=BRSK(U)$.
\end{eg}
\subsection{For any monomial $U$ in $\widetilde{\mathfrak{R}^v}\setminus \widetilde{\mathfrak{N}^v}$, $\tilde{\pi}(U)=BRSK(U)$}
Corollary~\ref{cor.a-main} below states that for any monomial $U$ in $\widetilde{\mathfrak{R}^v}\setminus \widetilde{\mathfrak{N}^v}$, $\tilde{\pi}(U)=BRSK(U)$. This is a corollary to Theorem~\ref{a.main}. But before we state Theorem~\ref{a.main} and Corollary~\ref{cor.a-main}, we need to define the map $\tilde{\pi}$ on monomials in $\widetilde{\mathfrak{R}^v}\setminus \widetilde{\mathfrak{N}^v}$. But for doing so, we need some preparation.\par
Let $v=(v_1,\ldots,v_d)$. Given $\beta_1=(r_1,c_1)$ and $\beta_2=(r_2,c_2)$ in$\widetilde{\mathfrak{R}^v}\setminus \widetilde{\mathfrak{N}^v}$, we say that $\beta_1>\beta_2$ if $r_1<r_2$ and $c_2<c_1$. A sequence $\beta_1>\cdots>\beta_t$ of elements of $\widetilde{\mathfrak{R}^v}\setminus \widetilde{\mathfrak{N}^v}$ is called an \textbf{anti}-$v$-\textbf{chain}. Given an anti-$v$-chain $\beta_1=(r_1,c_1)>\cdots>\beta_t=(r_t,c_t)$, we define
$$s_{\beta_1}\cdots s_{\beta_t}v:=(\{v_1,\ldots,v_d\}\setminus\{c_1,\ldots,c_t\})\cup\{r_1,\ldots,r_t\}.$$
We say that an element $w$ of $I(d,N)$ \textbf{anti-dominates} the anti-$v$-chain $\beta_1>\cdots>\beta_t$ if $w\leq s_{\beta_1}\cdots s_{\beta_t}v$. Let $\mathfrak{S}$ be a monomial in $\widetilde{\mathfrak{R}^v}\setminus \widetilde{\mathfrak{N}^v}$. We say that $w$ \textbf{anti-dominates} $\mathfrak{S}$ if $w$ anti-dominates every anti-$v$-chain in $\mathfrak{S}$.\par 
We call \textbf{distinguished} the subsets $\mathfrak{S}$ of $\widetilde{\mathfrak{R}^v}\setminus\widetilde{\mathfrak{N}^v}$ satisfying the following conditions:

\noindent (A) For $(r,c)\neq (r',c')$ in $\mathfrak{S}$, we have $r\neq r'$ and $c\neq c'$.\\
\noindent (B) If $\mathfrak{S}=\{(r_1,c_1),\ldots,(r_p,c_p)\}$ with $r_1>r_2>\ldots>r_p$, then for $j$, $1\leq j\leq p-1$, we have either $c_j<c_{j+1}$ or $r_j>c_{j+1}$.\\
\begin{remark}\label{r.anti-distinguished}
It can be proved similarly as in Proposition 4.3 of \cite{kr} that there exists a bijection between elements $w$ of $I(d,N)$ satisfying $w\leq v$ on the one hand and distinguished subsets of $\widetilde{\mathfrak{R}^v}\setminus\widetilde{\mathfrak{N}^v}$ on the other hand. We denote this bijective correspondence by $w\leftrightarrow\mathfrak{S}_w$. 
\end{remark}
Let $\mathfrak{S}$ be a non-empty monomial in $\widetilde{\mathfrak{R}^v}\setminus\widetilde{\mathfrak{N}^v}$. If $\beta_1>\cdots>\beta_t$ is an anti-$v$-chain in $\mathfrak{S}$, then we call $\beta_1$ the \textbf{head} of the anti-$v$-chain and $\beta_t$ its \textbf{tail}. We call $t$ to be the \textbf{length} of the anti-$v$-chain. We say that an element $\beta$ of $\mathfrak{S}$ is $t$-\textbf{deep} in $\mathfrak{S}$ (where $t$ is a positive integer) if $\beta$ is the tail of an anti-$v$-chain in $\mathfrak{S}$ of length $t$. The \textbf{depth} of $\beta$ in $\mathfrak{S}$ is defined to be $t$ if $\beta$ is $t$-deep in $\mathfrak{S}$ but not $(t+1)$-deep in $\mathfrak{S}$.\par
We will now define the map $\pi$ on any monomial in $\widetilde{\mathfrak{R}^v}\setminus\widetilde{\mathfrak{N}^v}$. Let $\mathfrak{S}$ be a non-empty monomial in the elements of $\widetilde{\mathfrak{R}^v}\setminus\widetilde{\mathfrak{N}^v}$. We partition $\mathfrak{S}$ in two stages. First we partition $\mathfrak{S}$ into subsets $\mathfrak{S}_1,\ldots,\mathfrak{S}_k$, where $k$ is the largest length of an anti-$v$-chain in $\mathfrak{S}$: $\beta\in\mathfrak{S}$ belongs to $\mathfrak{S}_j$ if it is $j$-deep but not $(j+1)$-deep.\par
Now we partition each $\mathfrak{S}_j$ into subsets called \textit{blocks} as follows. We arrange the elements of $\mathfrak{S}_j$ in non-increasing order of their row numbers (where elements occur with their respective multiplicities). Among those with the same row number, the arrangement is by non-increasing order of column numbers. Two consecutive members $(r,c)$, $(R,C)$ in this arrangement are said to be \textit{related} if $r<C$. The blocks are the equivalence classes of the smallest equivalence relation containing the above relations.\par
Let $\mathfrak{B}$ be a single block of some $\mathfrak{S}_j$. Let
$$(r_1,c_1),\ldots,(r_p,c_p)$$
be the elements of $\mathfrak{B}$ written in non-increasing order of both row and column numbers (in such an arrangement, the elements occur with their respective multiplicities). We set $w(\mathfrak{B}):=(r_p,c_1)$ and $\mathfrak{B}'$ to be the monomial
$$\{(r_1,c_2),(r_2,c_3),\ldots,(r_{p-2},c_{p-1}),(r_{p-1},c_p)\}.$$
Set $\mathfrak{S}_j^{(1)}:=\cup_{\mathfrak{B}}\mathfrak{B}'$ (where the index $\mathfrak{B}$ runs over all blocks of $\mathfrak{S}_j$) and $\mathfrak{S}^{(1)}:=\cup_{j=1}^{k}\mathfrak{S}_j^{(1)}$. It follows (similarly as in Corollary 4.13 of \cite{kr}) that the set
$$\{w(\mathfrak{B})|\mathfrak{B}\ \text{is a block of}\ \mathfrak{S}\}$$
is a distinguished subset of $\widetilde{\mathfrak{R}^v}\setminus\widetilde{\mathfrak{N}^v}$. Let $w$ be the corresponding element of $I(d,N)$ (under the correspondence given in Remark \ref{r.anti-distinguished}). Set 
$$\pi(\mathfrak{S}):=(w,\mathfrak{S}^{(1)}).$$
This finishes the description of the map $\pi$.\par 
A standard monomial $\theta_1\geq\cdots\geq\theta_t$ in $I(d,N)$ is called \textbf{dominated} by $v$ if $v\geq\theta_1$. Let $\widetilde{SM^v_v}$ denote the set of all $v$-compatible standard monomials in $I(d,N)$ dominated by $v$.\par
Using $\pi$, we now define the map $\tilde{\pi}$ from the set of all monomials in $\widetilde{\mathfrak{R}^v}\setminus\widetilde{\mathfrak{N}^v}$ to $\widetilde{SM^v_v}$. We proceed by induction on the degree of a monomial $\mathfrak{S}$ in $\widetilde{\mathfrak{R}^v}\setminus\widetilde{\mathfrak{N}^v}$. The image of the empty monomial under $\tilde{\pi}$ is taken to be the empty monomial. Let $\mathfrak{S}$ be non-empty, and suppose that $\pi(\mathfrak{S})=(w,\mathfrak{S}^{(1)})$. It can be shown (similarly as in (1) and (2) of Proposition 4.1 of \cite{kr}) that the degree of $\mathfrak{S}^{(1)}$ is strictly less than that of $\mathfrak{S}$, and so by induction $\tilde{\pi}(\mathfrak{S}^{(1)})$ is defined. Suppose that $\tilde{\pi}(\mathfrak{S}^{(1)})=w'\leq\cdots$. It can be shown (similarly as in (3) of Proposition 4.1 of \cite{kr}) that $w\leq w'$. We set $\tilde{\pi}(\mathfrak{S}):=w\leq\tilde{\pi}(\mathfrak{S}^{(1)})$. This finishes the description of the map $\tilde{\pi}$ on the set of all monomials in $\widetilde{\mathfrak{R}^v}\setminus\widetilde{\mathfrak{N}^v}$.\par 
The proof of Theorem \ref{a.main} below is similar to the proof of Theorem \ref{t.main}.
\begin{theorem}\label{a.main}
Let $U$ be a finite monomial in $\widetilde{\mathfrak{R}^v}\setminus \widetilde{\mathfrak{N}^v}$. Let $\pi(U)=(w_0,U^{(1)})$, $\pi(U^{(1)})=(w_1,U^{(2)})$, $\ldots$ and so on till $\pi(U^{(k)})=(w_k,\emptyset)$, where $\emptyset$ is the empty monomial. Then for each $r \in \{0,1, \ldots , k\}$, the following holds : \\
(i) All the row numbers of the distinguished subset $\mathfrak{S}_{w_r}$ (corresponding to $w_r$) consist of the $(r+1)$-th row entries of the left hand notched tableau of $BRSK(U)$.\\
(ii) All the column numbers of $\mathfrak{S}_{w_r}$ comprise of the $(r+1)$-th row entries of the right hand notched tableau of $BRSK(U)$.
\end{theorem}
Example~\ref{eg.tildepifornegative} below illustrates Theorem~\ref{a.main}.
\begin{eg}\label{eg.tildepifornegative}
Let $d=5,N=11$ and $v=(3,6,8,10,11)$. Let 
$$U=\{(9,11),(4,11),(7,10),(5,10),(7,8),(1,8),(4,6)\}$$
be a finite monomial in $\widetilde{\mathfrak{R}^v}\setminus\widetilde{\mathfrak{N}^v}$. 
\begin{figure}[h]
\setlength{\unitlength}{0.8mm}
\centering
\begin{picture}(60,70)(0,0)
\matrixput(0,0)(10,0){5}(0,10){1}{\circle*{.5}}
\matrixput(0,10)(10,0){4}(0,10){2}{\circle*{.5}}
\matrixput(0,30)(10,0){3}(0,10){1}{\circle*{.5}}
\matrixput(0,40)(10,0){2}(0,10){1}{\circle*{.5}}
\linethickness{.5pt}
\dottedline{0.1}(0,50)(0,20)(0,0)(20,0)
\dottedline{0.1}(10,40)(10,20)(30,20)
\put(0,20){\makebox(0,0){$\bullet$}}
\put(0,50){\makebox(0,0){$\bullet$}}
\put(10,30){\makebox(0,0){$\bullet$}}
\put(10,40){\makebox(0,0){$\bullet$}}
\put(20,0){\makebox(0,0){$\bullet$}}
\put(20,40){\makebox(0,0){$\bullet$}}
\put(30,20){\makebox(0,0){$\bullet$}}
\linethickness{0.3pt}
\multiput(0,0)(1.5,0){27}{\line(1,0){1}}
\multiput(0,10)(1.5,0){27}{\line(1,0){1}}
\multiput(0,20)(1.5,0){20}{\line(1,0){1}}
\multiput(0,30)(1.5,0){20}{\line(1,0){1}}
\multiput(0,40)(1.5,0){14}{\line(1,0){1}}
\multiput(0,50)(1.5,0){7}{\line(1,0){1}}
\multiput(0,0)(0,1.5){34}{\line(0,1){1}}
\multiput(10,0)(0,1.5){34}{\line(0,1){1}}
\multiput(20,0)(0,1.5){27}{\line(0,1){1}}
\multiput(30,0)(0,1.5){20}{\line(0,1){1}}
\multiput(40,0)(0,1.5){7}{\line(0,1){1}}
\multiputlist(-10,0)(0,10){1,2,4,5,7,9}
\put(-3,54){11}
\put(8,54){10}
\put(20,43){8}
\put(30,33){6}
\put(40,13){3}
\put(2,15){1}
\put(2,45){1}
\put(12,25){1}
\put(6,35){1}
\put(23,2){1}
\put(15,35){1}
\put(25,15){1}
\end{picture}
\caption{The monomial $U$ in $\widetilde{\mathfrak{R}^v}\setminus\widetilde{\mathfrak{N}^v}$ and its block decomposition}
\label{figure.6}
\end{figure}
Figure~\ref{figure.6} shows the monomial $U$ and its block decomposition. The dark circles in the figure represent points in the monomial $U$ with their respective multiplicities (the multiplicity of each point in the monomial $U$ is $1$ here, which is written near those points in the grid). The dark line segments (together with the point $(7,8)$) denote the blocks of $U$.\par
For this monomial $U$, we have $\pi(U)=(w_0,U^{(1)})$, where $w_0=(1,3,4,6,7)$ and $U^{(1)}=\{(9,11),(4,8),(7,10),(5,6)\}$. Then $\pi(U^{(1)})=(w_1,U^{(2)})$, where $w_1=(3,4,5,8,10)$ and $U^{(2)}=\{(9,10),(7,8)\}$. And finally, $\pi(U^{(2)})=(w_2,\emptyset)$, where $w_2=(3,6,7,9,11)$.\par 
So we have $k=2$, $\mathfrak{S}_{w_0}=\{(1,11),(4,10),(7,8)\}$, $\mathfrak{S}_{w_1}=\{(4,11),(5,6)\}$ and $\mathfrak{S}_{w_2}=\{(9,10),(7,8)\}$. It is also easy to check that
\newcommand{\J}{10}
\newcommand{\C}{11}
$$BRSK(U)=\left(\hspace{0.2cm}\young(147,45,79),\hspace{0.5cm}\young(8\J\C,6\C,8\J)\hspace{0.2cm}\right).$$
Note that for $r=0$, all the row numbers of $\mathfrak{S}_{w_0}$ are given by $\{1,4,7\}$, which consist of the first row entries of the left hand notched tableau of $BRSK(U)$. And all the column numbers of $\mathfrak{S}_{w_0}$ are given by $\{8,10,11\}$, which consist of the first row entries of the right hand notched tableau of $BRSK(U)$. Similarly for $r=1$ and $r=2$.
\end{eg}
\begin{corollary}\label{cor.a-main}
For any monomial $U$ in $\widetilde{\mathfrak{R}^v}\setminus \widetilde{\mathfrak{N}^v}$, $\tilde{\pi}(U)=BRSK(U)$.
\end{corollary}
\begin{corollary}
So by Corollary \ref{t.main} and Corollary \ref{a.main}, we can say that for any monomial $U$ in $\widetilde{\mathfrak{R}^v}$, $\tilde{\pi}(U)=BRSK(U)$.
\end{corollary}
\subsection{Some results needed to prove theorem \ref{t.main}}\label{ss.results}
\noindent We will prove Theorem \ref{t.main} by induction on the cardinality $n$ of the monomial $U$. If $n=1$, then the theorem is obvious. To prove Theorem \ref{t.main} above, we need some notation, lemmas and definitions. We will mention those first. In the rest of this paper, we are going to assume the  \textbf{induction hypothesis}, that is, Theorem \ref{t.main} holds true for all finite monomials in $\widetilde{\mathfrak{N}^v}$ of cardinality $\leq n-1$. Also, in the rest of this paper, (unless otherwise mentioned) we are going to use the same terminology and notation as in the papers \cite{kr} and \cite{kv}.\par
\noindent\textbf{Notation:} Let $\iota$ be the involution map, which was defined in \S\ref{ss.recap}. Let $U$ be a monomial in $\widetilde{\mathfrak{N}^v}$ of degree $n$. Arrange $ \iota(U)$ in lexicographic order, say,
$\iota(U)=\{(a_1,b_1),(a_2,b_2),\ldots,(a_n,b_n)\}$.
Let $F=\{(b_1,a_1),(b_2,a_2),\ldots,(b_{n-1},a_{n-1})\}$. Then $b_n \leq b_{n-1}$ and if $b_n = b_{n-1}$, we have $a_{n-1} \geq a_n$. The element $(b_n,a_n)$ enters into $F$ to make it $U$. Let $BRSK(\iota(F))=(P^{(n-1)},Q^{(n-1)})$ .\\
Let $p_{ij}$ denote the entry in the $i$-th row and $j$-th column of $P^{(n-1)}$. Similarly, let $q_{ij}$ denote the entry in the $i$-th row and $j$-th column of $Q^{(n-1)}$. Clearly, $BRSK(\iota(U))= (P^{(n-1)},Q^{(n-1)})\la{b_n}a_n$, the entries of $P^{(n-1)}$  are $\{a_1,\ldots,a_{n-1}\}$ and the entries of $Q^{(n-1)}$ are $\{b_1,\ldots,b_{n-1}\}$. Let $(P^{(n)},Q^{(n)})$ denote $(P^{(n-1)},Q^{(n-1)})\la{b_n}a_n$. Let $k$ be a positive integer such that elements of depth $k$ exist in the monomial $F$. Let $\{(r_1,c_1),\ldots,(r_p,c_p)\}$ denote the topmost block of $F$ of depth $k$, where the elements of the block are written in non-decreasing order of both row and column numbers. It then follows from the induction hypothesis that $c_1$ is an entry in the first row of $P^{(n-1)}$.
\begin{lemma}\label{l.A}
The entries in the first row of $P^{(n-1)}$ which are strictly less than $c_1$, are the column numbers of the topmost elements of some blocks of $F$ of depth $<k$. 
\end{lemma}
\begin{Proof}
Suppose not. Say, some $p_{1j^\prime} (< c_1)$ is the column number of the topmost element of some block $\mathfrak{B}$ of $F$ of depth ${\geq k}$. Say, the first element of the block $\mathfrak{B}$ is $(R,p_{1j^\prime})$.\par
\noindent \textbf{Case (i) : $R> r_1$} \\
In this case, $(r_1,c_1)$ will have depth $> k$ in $F$, a contradiction.\\
\textbf{Case (ii) : $R=r_1$}\\
In this case, $(r_1,c_1)$ will either have depth $> k$ (a contradiction) or both $(r_1,c_1)$ and $(R,p_{1j^\prime})$ have depth equal to $k$ in $F$. But in the later situation, $(r_1,c_1)$ cannot be the first element in the topmost block of $F$ of depth $k$, again a contradiction.\\
\textbf{Case (iii) : $R< r_1$}\\
If $(R,p_{1,j^\prime})$ is of depth $k$, then clearly $(r_1,c_1)$ cannot be the first element of the topmost block of $F$ of depth $k$. If $(R,p_{1j^\prime})$ has depth $s >k$, then $\exists$ a $v$-chain of length $s$ having tail $(R,p_{1j^\prime})$. Say, the $v$-chain is $(e_1,f_1)>(e_2,f_2)> \cdots >(e_s,f_s)=(R,p_{1j^\prime})$. Then $(e_k,f_k)$ will have depth $k$ in $F$. If $e_k \leq r_1$, then actually $(r_1,c_1)$ is not the topmost element of the topmost block of $F$ of depth $k$, a contradiction. If $e_k > r_1$, then $(e_k,f_k)>(r_1,c_1)$, a contradiction to the depth of $(r_1,c_1)$ in $F$.
\end{Proof}
Example~\ref{eg.lemmaA} below illustrates Lemma~\ref{l.A}.
\begin{eg}\label{eg.lemmaA}
Let $d=9$, $N=17$, and $v=(1,2,4,6,7,11,12,13,15)$. Let $F$ be the finite monomial in $\widetilde{\mathfrak{N}^v}$ given by
$$F=\{(8,2),(8,2),(8,6),(8,6),(8,7),(9,2),(10,1),(14,13),(16,12),(16,13),(17,11)\}.$$ 
Figure~\ref{figure.4} shows the monomial $F$ and its block decomposition. The dark circles denote the elements of the monomial $F$, and the numbers written near these dark circles denote the multiplicities with which these elements occur in the monomial $F$. The dark line segments (in case of blocks consisting of more than one distinct elements) and the dark circles (in case of blocks consisting of a single distinct element) are the blocks of $F$. 
\begin{figure}[h]
\setlength{\unitlength}{0.6mm}
\centering
\begin{picture}(60,80)(0,0)
\matrixput(0,0)(10,0){7}(0,10){3}{\circle*{.5}}
\matrixput(0,30)(10,0){4}(0,10){2}{\circle*{.5}}
\matrixput(0,40)(10,0){2}(0,10){2}{\circle*{.5}}
\linethickness{.7pt}
\dottedline{0.1}(10,50)(10,40)
\linethickness{.9pt}
\dottedline{0.1}(60,10)(70,10)
\dottedline{0.1}(30,50)(40,50)
\put(0,30){\makebox(0,0){$\bullet$}}
\put(10,40){\makebox(0,0){$\bullet$}}
\put(10,50){\makebox(0,0){$\bullet$}}
\put(30,50){\makebox(0,0){$\bullet$}}
\put(40,50){\makebox(0,0){$\bullet$}}
\put(50,0){\makebox(0,0){$\bullet$}}
\put(60,10){\makebox(0,0){$\bullet$}}
\put(70,10){\makebox(0,0){$\bullet$}}
\put(70,20){\makebox(0,0){$\bullet$}}
\linethickness{0.1pt}
\multiput(0,0)(1.5,0){54}{\line(1,0){1}}
\multiput(0,10)(1.5,0){54}{\line(1,0){1}}
\multiput(0,20)(1.5,0){47}{\line(1,0){1}}
\multiput(0,30)(1.5,0){27}{\line(1,0){1}}
\multiput(0,40)(1.5,0){27}{\line(1,0){1}}
\multiput(0,50)(1.5,0){27}{\line(1,0){1}}
\multiput(0,60)(1.5,0){14}{\line(1,0){1}}
\multiput(0,70)(1.5,0){7}{\line(1,0){1}}
\multiput(0,0)(0,1.5){47}{\line(0,1){1}}
\multiput(10,0)(0,1.5){47}{\line(0,1){1}}
\multiput(20,0)(0,1.5){40}{\line(0,1){1}}
\multiput(30,0)(0,1.5){34}{\line(0,1){1}}
\multiput(40,0)(0,1.5){34}{\line(0,1){1}}
\multiput(50,0)(0,1.5){14}{\line(0,1){1}}
\multiput(60,0)(0,1.5){14}{\line(0,1){1}}
\multiput(70,0)(0,1.5){14}{\line(0,1){1}}
\multiput(80,0)(0,1.5){7}{\line(0,1){1}}
\multiputlist(-10,0)(0,10){17,16,14,10,9,8,5,3}
\put(0,75){1}
\put(10,75){2}
\put(20,65){4}
\put(30,55){6}
\put(40,55){7}
\put(48,25){11}
\put(58,25){12}
\put(68,25){13}
\put(78,15){15}
\put(2,32){1}
\put(12,35){1}
\put(12,52){2}
\put(26,45){2}
\put(36,45){1}
\put(50,2){1}
\put(56,11){1}
\put(71,5){1}
\put(66,15){1}
\end{picture}
\caption{The monomial $F$ and its block decomposition}
\label{figure.4}
\end{figure}
Let $BRSK(\iota(F))=(P^{(n-1)},Q^{(n-1)})$.\\ Then 
\newcommand{\J}{10}
\newcommand{\B}{12}
\newcommand{\A}{13}
\newcommand{\C}{11}
\newcommand{\D}{14}
\newcommand{\E}{16}
\newcommand{\F}{17}
$$(P^{(n-1)},Q^{(n-1)})=\left(\hspace{0.2cm}\young(126\C\B\A,2\A,2,6,7),\hspace{0.5cm}\young(89\J\D\E\F,8\E,8,8,8)\hspace{0.2cm}\right).$$\\ \\
Let $k=3$. Then $\{(r_1,c_1),(r_2,c_2),(r_3,c_3)\}$ is the topmost block of $F$ of depth $k$, where $(r_1,c_1)=(r_2,c_2)=(8,6)$ and $(r_3,c_3)=(8,7)$. Here $c_1=6$. Observe that $c_1=6$ is an entry in the first row of $P^{(n-1)}$. Observe also that the entries in the first row of $P^{(n-1)}$ which are strictly less than $c_1$ are $1$ and $2$. Note that $1$ and $2$ are the column numbers of the topmost elements of some blocks of $F$ of depth $<k=3$. In fact, $1$ is the column number of the topmost element of the block $\{(10,1)\}$ of $F$, which is of depth $1$. Similarly, $2$ is the column number of the topmost element of the block $\{(8,2),(8,2),(9,2)\}$ of $F$, which is of depth $2$.  
\end{eg}
\begin{lemma}\label{l.B}
Let $(b_n,a_n)$ enter into the monomial $F$ to make it $U$ in such a way that the singleton set $\{(b_n,a_n)\}$ is the topmost block of $U$ of depth $k$, and the next block of $U$ of depth $k$ from the top being $\{(r_1,c_1),(r_2,c_2),\ldots,(r_p,c_p)\}$. Then the entries in the first row of $P^{(n-1)}$ which are strictly less than $b_n$ are all strictly less than $a_n$.
\end{lemma}
\begin{Proof}
It follows from the hypothesis of this lemma that $b_n<c_1$, $b_n<r_1\leq r_2\leq\cdots\leq r_p$ and $a_n<c_1\leq c_2\leq\cdots\leq c_p$. It also follows from the induction hypothesis of Theorem \ref{t.main} that the first row of $P^{(n-1)}$ contains the smallest column numbers of each block of $F$. In particular, it contains the entry $c_1$.\par
We will prove this lemma by method of contradiction. Suppose the conclusion of this lemma does not hold. Say, some entry $p_{1j_0}$ of the first row of $P^{(n-1)}$ which is strictly less than $b_n$ is $\geq a_n$. Then we have $a_n \leq p_{1j_0} <b_n < c_1$. By induction hypothesis we know that $p_{1j_0}$ is the smallest column number of some block of $F$, say block $\mathfrak{D}$.
Since $P_{1j_0}< c_1$, it follows (from Lemma \ref{l.A}) that the block $\mathfrak{D}$ has depth $<k$.
 Say, $\mathfrak{D}$ is a block of depth $s(<k)$. Now, since $(b_n,a_n)$ is of depth $k$ in $U$, there exists an element $(R,C)$ in $U$ of depth $s\ (<k)$ such that $(R,C)$ and $(b_n,a_n)$ form a $v$-chain. That is, $b_n < R$ and $C<a_n$. Say, $(R,C)$ lies in the block $\mathfrak{B}$ of $U$. Clearly then, $\mathfrak{B} \neq \mathfrak{D}$ (because the smallest column number of the block $\mathfrak{D}$ is $P_{ij_0}$ which is $\geq a_n$). Let $(\widehat{R}, \widehat{C})$ be the bottom-most element of the block $\mathfrak{B}$. Then $R \leq \widehat{R}$ and $C\leq \widehat{C}$.\par
Hence, we have $\widehat{R}\geq R>b_n>P_{ij_0} \Rightarrow \widehat{R} >P_{ij_0}$. That is, $\mathfrak{B}$ and $\mathfrak{D}$ are not two different blocks of depth $s$, a contradiction. 
\end{Proof}
Example~\ref{eg.lemmaB} below illustrates Lemma~\ref{l.B}.
\begin{eg}\label{eg.lemmaB}
Let $d,N,v,F$ be as given in Example~\ref{eg.lemmaA} above. Then $(P^{(n-1)},Q^{(n-1)})$ is also given in Example~\ref{eg.lemmaA} above. Let $(b_n,a_n)=(5,4)$ and $k=3$. Then $(b_n,a_n)$ enters into the monomial $F$ to make it $U$ in such a way that the singleton set $\{(b_n,a_n)\}$ is the topmost block of $U$ of depth $k$. Observe that the entries in the first row of $P^{(n-1)}$ which are strictly less than $b_n=5$ are $1$ and $2$. Both $1$ and $2$ are strictly less than $a_n=4$ as well.
\end{eg}
Let $\{(R_1, C_1),(R_2,C_2),\ldots,(R_p,C_p)\}$ be a block $\mathfrak{B}$ of some finite monomial $\mathfrak{S}$ of $\widetilde{\mathfrak{N}^v}$. Let $b_0 \leq R_1$ be such that $b_0 \in [N] \setminus v$ and $b_0 > C_1$. Let $a_0 \in v$ be such that $a_0 \leq C_1$. Then we say that $\{(b_0,a_0),(R_1,C_1),....,(R_p,C_p)\}$ is a \textbf{left concatenation} of $\mathfrak{B}$ by $(b_0,a_0)$.\par 
Example~\ref{eg.leftconcatenation} below illustrates the above definition of left concatenation of a block.
\begin{eg}\label{eg.leftconcatenation}
Let $d,N,v$ be as in Example~\ref{eg.lemmaA} above. Let $\mathfrak{S}$ be the monomial $F$ of $\widetilde{\mathfrak{N}^v}$ as given in Example~\ref{eg.lemmaA} above. Consider the block $\mathfrak{B}$ of $\mathfrak{S}$ given by $\mathfrak{B}=\{(8,6),(8,6),(8,7)\}$. Then $p=3$, $(R_1,C_1)=(R_2,C_2)=(8,6)$ and $(R_3,C_3)=(8,7)$. Let $(b_0,a_0)=(8,4)$. Then $b_0\leq R_1$, $b_0\in [N]\setminus v$ and $b_0>C_1$. Also, $a_0\in v$ and $a_0\leq C_1$. Observe that $\{(8,4),(8,6),(8,6),(8,7)\}$ is a left concatenation of the block $\mathfrak{B}$ by $(b_0,a_0)$. In Figure~\ref{figure.4} (inside Example~\ref{eg.lemmaA}), there was the block $\mathfrak{B}=\{(8,6),(8,6),(8,7)\}$. In Figure~\ref{figure.5}, we can see that $\{(8,4),(8,6),(8,6),(8,7)\}$ is a left concatenation of the block $\mathfrak{B}$ by $(8,4)$. 
\begin{figure}[h]
\setlength{\unitlength}{0.6mm}
\centering
\begin{picture}(60,110)(0,0)
\matrixput(0,0)(10,0){7}(0,10){3}{\circle*{.5}}
\matrixput(0,30)(10,0){4}(0,10){2}{\circle*{.5}}
\matrixput(0,40)(10,0){2}(0,10){2}{\circle*{.5}}
\linethickness{.7pt}
\dottedline{0.1}(10,50)(10,40)
\linethickness{.9pt}
\dottedline{0.1}(60,10)(70,10)
\dottedline{0.1}(20,50)(40,50)
\put(0,30){\makebox(0,0){$\bullet$}}
\put(10,40){\makebox(0,0){$\bullet$}}
\put(10,50){\makebox(0,0){$\bullet$}}
\put(30,50){\makebox(0,0){$\bullet$}}
\put(40,50){\makebox(0,0){$\bullet$}}
\put(50,0){\makebox(0,0){$\bullet$}}
\put(60,10){\makebox(0,0){$\bullet$}}
\put(70,10){\makebox(0,0){$\bullet$}}
\put(70,20){\makebox(0,0){$\bullet$}}
\put(20,50){\makebox(0,0){$\bullet$}}
\linethickness{0.1pt}
\multiput(0,0)(1.5,0){54}{\line(1,0){1}}
\multiput(0,10)(1.5,0){54}{\line(1,0){1}}
\multiput(0,20)(1.5,0){47}{\line(1,0){1}}
\multiput(0,30)(1.5,0){27}{\line(1,0){1}}
\multiput(0,40)(1.5,0){27}{\line(1,0){1}}
\multiput(0,50)(1.5,0){27}{\line(1,0){1}}
\multiput(0,60)(1.5,0){14}{\line(1,0){1}}
\multiput(0,70)(1.5,0){7}{\line(1,0){1}}
\multiput(0,0)(0,1.5){47}{\line(0,1){1}}
\multiput(10,0)(0,1.5){47}{\line(0,1){1}}
\multiput(20,0)(0,1.5){40}{\line(0,1){1}}
\multiput(30,0)(0,1.5){34}{\line(0,1){1}}
\multiput(40,0)(0,1.5){34}{\line(0,1){1}}
\multiput(50,0)(0,1.5){14}{\line(0,1){1}}
\multiput(60,0)(0,1.5){14}{\line(0,1){1}}
\multiput(70,0)(0,1.5){14}{\line(0,1){1}}
\multiput(80,0)(0,1.5){7}{\line(0,1){1}}
\multiputlist(-10,0)(0,10){17,16,14,10,9,8,5,3}
\put(0,75){1}
\put(10,75){2}
\put(20,65){4}
\put(30,55){6}
\put(40,55){7}
\put(48,25){11}
\put(58,25){12}
\put(68,25){13}
\put(78,15){15}
\put(16,55){1}
\put(2,32){1}
\put(12,35){1}
\put(5,51){2}
\put(26,45){2}
\put(36,45){1}
\put(50,2){1}
\put(56,11){1}
\put(71,5){1}
\put(66,15){1}
\end{picture}
\caption{Left concatenation of a block}
\label{figure.5}
\end{figure}
\end{eg}
\begin{lemma}\label{l.C}
Let $\{(r_1,c_1),\ldots,(r_p,c_p)\}$ be the topmost block of $F$ of depth $k$. Let $(b_n,a_n)$ enter into the monomial $F$ to make it $U$ in such a way that $\{(b_n,a_n),(r_1,c_1),\ldots,(r_p,c_p)\}$ becomes the topmost block of $U$ of depth $k$. Let $m$ be the positive integer as given in the statement of Theorem \ref{t.main}. Let $U^{(0)}:=U$. Then $\exists$ an integer $k'$ where $0\leq k' \leq m-1$ such that:\\
\textbf{(i)} For each $t\in\{0,1,\dots,k'\}$, all the blocks of $U^{(t)}$ except one are the same as the blocks of $F^{(t)}$. The one block of $U^{(t)}$ that is different, is, in fact, a left concatenation of a block of $F^{(t)}$ by $(b_n,\star)$, where $\star$ is some entry of $v$ which is $\geq a_n$.\\
\noindent \textbf{(ii)} The set of all blocks of $U^{(k'+1)}$ is equal to the set of all blocks of $F^{(k'+1)}$ union one more block, which is of the form $\{(b_n,\star)\}$, where $\star$ is some entry of $v$ which is $\geq a_n$.\\
\noindent\textbf{(iii)} For each $t \in \{k'+2,\ldots,m\}$, the set of all blocks of $U^{(t)}$ is the same as the set of all blocks of $F^{(t)}$.
\end{lemma}
\begin{Proof}
Clearly all blocks of $U$ except one (namely, $\{(b_n,a_n),(r_1,c_1),\ldots,(r_p,c_p)\}$) are the same as all blocks of $F$. And the block  $\{(b_n,a_n),(r_1,c_1),\ldots,(r_p,c_p)\}$ is a left concatenation of the block  $\{(r_1,c_1),\ldots,(r_p,c_p)\}$ of $F$ by $(b_n,a_n)$.\\
Let $\mathfrak{B}_1,\mathfrak{B}_2,\ldots,\mathfrak{B}_{n_0}$ denote the all other blocks of $U$ or $F$ (they are the same!). $U^{(1)}$ contains the elements  $(b_n,c_1),(r_1,c_2),\ldots,(r_{p-1},c_p)$, and all other elements of $U^{(1)}$ are given by $\mathfrak{B}_1^\prime,\ldots,\mathfrak{B}_{n_0}^\prime$. If $\{(b_n,c_1)\}$ is a single block in $U^{(1)}$, then $k'=0$.\par 
Now suppose $\{(b_n,c_1)\}$ is not a single block in $U^{(1)}$.

\noindent\textbf{Claim:} There exists a block in $U^{(1)}$, which is a left concatenation of a block of $F^{(1)}$ by $(b_n,c_1)$.

\noindent\textbf{Proof of claim:} Suppose not. Then $(b_n,c_1)$ is not the topmost element of any block of $U^{(1)}$. Since $(b_n,c_1)\in U^{(1)}$, there exists some block (say, $\mathfrak{D}$) of $U^{(1)}$ such that $(b_n,c_1)\in\mathfrak{D}$, $\mathfrak{D}=\{(\widehat{r_1},\widehat{c_1}),\ldots,(\widehat{r_k},\widehat{c_k})\}$, and $(b_n,c_1)=(\widehat{r_j},\widehat{c_j})$ for some $j\in\{2,\ldots,k\}$. Since $b_n\leq b_i$ for all $i\in\{1,\ldots,n-1\}$, we have $\widehat{r_1}=\cdots=\widehat{r_j}=b_n$. Also since $\{(\widehat{r_1},\widehat{c_1}),\ldots,(\widehat{r_k},\widehat{c_k})\}$ is a block (of $U^{(1)}$), we have $\widehat{c_1}\leq\cdots\leq\widehat{c_j}=c_1$. Since $\widehat{r_1}=\cdots=\widehat{r_j}=b_n$ and $(b_n,c_1)$ is not the topmost element of any block of $U^{(1)}$, we must have $\widehat{c_1}<c_1$. Observe that the element $(\widehat{r_1},\widehat{c_1})$ actually belongs to $F^{(1)}$ (because $U^{(1)}$ and $F^{(1)}$ differ only by one element, namely $(b_n,c_1)$ and $\widehat{c_1}<c_1$). This implies that there had existed an element $(b_n,\widehat{p})$ in $F$ such that $\widehat{p}\leq\widehat{c_1}$. Clearly $\widehat{p}\geq a_n$ (because $a_n$ is the least column number of any element in $U$ having row number $b_n$). So we have
$$a_n\leq\widehat{p}\leq\widehat{c_1}<c_1.$$
In particular, we have $a_n<c_1$.\\
\noindent Case (I): If $\widehat{p}=a_n$.\\
\noindent In this case, we have $(b_n,\widehat{p})=(b_n,a_n)\in F$. Clearly then $b_n=b_{n-1}$ and $a_n=a_{n-1}$. Since $b_n\leq b_i$ for all $i\in\{1,\ldots,n-1\}$ and $a_{n-1}$ ($=a_n$) is the least column number of any element in $F$ having row number $b_{n-1}$ ($=b_n$), we must have that $(b_n,\widehat{p})=(b_n,a_n)$ is the topmost element of the topmost block of $F$ of some depth. So when $(b_n,a_n)$ enters into the monomial $F$ to make it $U$, it gets added to the same block of $F$ to which $(b_n,\widehat{p})=(b_n,a_n)$ belongs. Therefore the block $\{(r_1,c_1),\ldots,(r_p,c_p)\}$ of $F$, to which $(b_n,a_n)$ gets added, must have the property that $(r_1,c_1)=(b_n,a_n)$. This implies that $a_n=c_1$, which contradicts the fact that $a_n<c_1$.\\
\noindent Case (II): If $\widehat{p}>a_n$.\\
\noindent In this case, we have $a_n<\widehat{p}<c_1$. Since $\{(b_n,a_n),(r_1,c_1),\ldots,(r_p,c_p)\}$ is a block of $U$, we have $b_n\leq r_1$ and $b_n>c_1$. Now since $(b_n,\widehat{p})\in F$, $b_n\leq r_1$, $\widehat{p}<c_1$ and $b_n>c_1$, we must have that either $(b_n,\widehat{p})$ belongs to the same block (of depth $k$) of $F$ as $(r_1,c_1)$ or the depth of $(b_n,\widehat{p})$ in $F$ is not equal to $k$.\par 
If $(b_n,\widehat{p})$ belongs to the same block (of depth $k$) of $F$ as $(r_1,c_1)$, then since $b_n\leq r_1$ and $\widehat{p}<c_1$, the element $(b_n,\widehat{p})$ must come before the element $(r_1,c_1)$ in that block of $F$. But this is a contradiction because we already know that $(r_1,c_1)$ is the topmost element of some block of $F$.\par
If the depth of $(b_n,\widehat{p})$ in $F$ is $<k$, then let $s$ ($<k$) be the depth of $(b_n,\widehat{p})$ in $F$. Since $(b_n,a_n)$ is the topmost element of the topmost block of $U$ of depth $k$, there must exist an element $(g,h)\in F$ such that $(g,h)>(b_n,a_n)$ and the depth of $(g,h)$ in $F$ is $s$. Then $g>b_n$ and $h<a_n$. This implies that $g>b_n$ and $h<\widehat{p}$, which in turn implies that $(g,h)>(b_n,\widehat{p})$. This is a contradiction because two elements in $F$ of the same depth cannot be comparable.\par
If the depth of $(b_n,\widehat{p})$ in $F$ is $>k$, then there must exist an element $(e,f)$ in $F$ of depth $k$ such that $(e,f)>(b_n,\widehat{p})$. So we have $e>b_n$ and $f<\widehat{p}$. Since $f<\widehat{p}$ and $\widehat{p}<c_1$, we have $(e,f)\neq (r_1,c_1)$. Also $b_n>c_1>\widehat{p}>f$. So $(b_n,a_n)$ and $(e,f)$ belong to the same block of $U$ of depth $k$ and $(b_n,a_n)$ comes before $(e,f)$ in that block of $U$. Note that if $e>r_1$, then we will have $(e,f)>(r_1,c_1)$, which is impossible. Hence $e\leq r_1$. Since $e\leq r_1$ and $f<c_1$, the element $(e,f)$ comes before the element $(r_1,c_1)$ in the block of $U$ of depth $k$ containing $(r_1,c_1)$. The above implies that $(e,f)$ lies strictly between $(b_n,a_n)$ and $(r_1,c_1)$ in this block of $U$, which is not possible. So we arrive at a contradiction. This proves our claim.\par
So, the blocks of $U^{(1)}$ and $F^{(1)}$ are essentially the same except one block, which is a left concatenation of a block of  $F^{(1)}$ by $(b_n,c_1)$. Let that block of $U^{(1)}$ be $\{(b_n,c_1),(\widetilde{r_1},\widetilde{c_2}),\ldots,(\widetilde{r_{l-1}},\widetilde{c_l})\}=\mathfrak{C}$. Then $\mathfrak{C}^\prime=\{(b_n,\widetilde{c_2}),(\widetilde{r_1},\widetilde{c_3}),\ldots,(\widetilde{r_{l-2}},\widetilde{c_l})\}$, where $c_1\leq \widetilde{c_2}$. Now, $U^{(2)}$ contains the elements of $\mathfrak{C}^\prime$. If $\{(b_n,\widetilde{c_2})\}$ is a single block in $U^{(2)}$, take $k'=1$. Otherwise proceed similarly. This process will stop at a stage $k'\leq m-1$, because $U$ is a finite monomial, and at some stage, we will surely get a block consisting of a single element $\{(b_n,\star)\}$, where $\star$ is some entry of $v$ which is $\geq a_n$. And after the $(k'+1$)-th stage, again the blocks of $U^{(t)}$ and $F^{(t)}$ will remain the same.
\end{Proof}
\begin{remark}
 Lemma \ref{l.C} above simply means that in the process \\$(P^{(n-1)},Q^{(n-1)})\la{b_n}a_n$, there is bumping in $P^{(n-1)}$ upto the $k'$-th stage. At the $(k'+1)$-th stage, the bumping stops and $b_n$ is placed in a new box in some row of $Q^{(n-1)}$. The remaining entries of $P^{(n-1)}$ and $Q^{(n-1)}$ will remain unchanged.
\end{remark}
Example~\ref{eg.lemmaC} below illustrates Lemma~\ref{l.C}.
\begin{eg}\label{eg.lemmaC}
Let $d,N,v,F$ be as in Example~\ref{eg.lemmaA} above. Let $k=2$. Then $\{(8,2),(8,2),(9,2)\}$ is the topmost block of $F$ of depth $k$. Let $(b_n,a_n)=(5,2)$. Then $(b_n,a_n)=(5,2)$ enters into the monomial $F$ to make it $U$ in such a way that $\{(5,2),(8,2),(8,2),(9,2)\}$ becomes the topmost block of $U$ of depth $k$. Observe that $\pi(U)=(w_0,U^{(1)}),\pi(U^{(1)})=(w_1,U^{(2)}),\pi(U^{(2)})=(w_2,U^{(3)}),\pi(U^{(3)})=(w_3,U^{(4)})$ and $\pi(U^{(4)})=(w_4,\emptyset)$, where $w_0,w_1,w_2,w_3,w_4$ are some elements of $I(9,17)$, $U^{(1)}=\{(5,2),(8,2),(8,2),(16,13),(8,6),(8,7)\}$, $U^{(2)}=\{(5,2),(8,2),(8,6),(8,7)\}$, $U^{(3)}=\{(5,2),(8,6),(8,7)\}$, and $U^{(4)}=\{(8,7)\}$. So $m=4$ here.\par 
Observe that there exists an integer $k'=2\in\{0,1,2,3(=m-1)\}$ such that the following holds true:\\
\noindent (i) All the blocks of $U^{(0)}(=U)$ except one are the same as the blocks of $F^{(0)}(=F)$. The one block of $U^{(0)}$ that is different, is, in fact, a left concatenation of the block $\{(8,2),(8,2),(9,2)\}$ of $F^{(0)}$ by $(5,2)$. Similarly, all the blocks of $U^{(1)}$ except one are the same as the blocks of $F^{(1)}$. The one block of $U^{(1)}$ that is different, is, in fact, a left concatenation of the block $\{(8,2),(8,2),(8,6),(8,7)\}$ of $F^{(1)}$ by $(5,2)$. Finally, all the blocks of $U^{(2)}$ except one are the same as the blocks of $F^{(2)}$. The one block of $U^{(2)}$ that is different, is, in fact, a left concatenation of the block $\{(8,2),(8,6),(8,7)\}$ of $F^{(2)}$ by $(5,2)$.\\
\noindent (ii) The set of all blocks of $U^{(3)}$ is equal to the set of all blocks of $F^{(3)}$ plus one more block, which is $\{(5,2)\}$.\\
\noindent (iii) The set of all blocks of $U^{(4)}$ is the same as the set of all blocks of $F^{(4)}$.
\end{eg}
\begin{lemma}\label{l.D}
Let $\{(r_1,c_1),\ldots,(r_p,c_p)\}$ be the topmost block of $F$ of depth $k$. Let $(b_n,a_n)$ enter into the monomial $F$ to make it $U$ in such a way that $\{(b_n,a_n),(r_1,c_1),\ldots,(r_p,c_p)\}$ becomes the topmost block of $U$ of depth $k$. Let $(P^{(n)},Q^{(n)})$ = $BRSK(\iota(U))=(P^{(n-1)},Q^{(n-1)})\la{b_n}a_n$. Then \\
\noindent \textbf{Step I:} $a_n$ is an entry in the first row of $P^{(n)}$ and $r_p$ is an entry of the first row of $Q^{(n)}$. $a_n$ bumps $c_1$ from the first row of $P^{(n-1)}$, $c_1$ is placed in the second row of $P^{(n-1)}$. \\
\noindent \textbf{Step II:} $c_1$ either bumps an entry from the 2nd row of $P^{(n-1)}$ (Case I) or it does not bump anything (Case II). The bumping by $c_1$ happens iff $\exists$ a block in $U^{(1)}$ having at least two elements and having $c_1$ as the least column number.\par
This process of bumping continues upto a finite stage until $\{(b_n, \star)\}$ is a single block, where $\star$ is some entry of $v$ which is $\geq a_n$. At this stage, where $\{(b_n,\star)\}$ is a single block, $\star$ is placed in a new box in some row of $P^{(n-1)}$ and $b_n$ is placed in a new box in the corresponding row of $Q^{(n-1)}$. All other entries of $P^{(n-1)}$ and $Q^{(n-1)}$ remain unchanged in the same rows.
\end{lemma}
Before going to the proof of Lemma \ref{l.D}, let us have a pictorial look of it and an example which illustrates it. Figure~\ref{figure.7} below provides a pictorial illustration of Lemma~\ref{l.D} and Example~\ref{eg.lemmaD} illustrates it.
\begin{figure}
\setlength{\unitlength}{0.6mm}
\begin{picture}(60,100)(0,0)
\put(-15,60){{Step I:}}
\put(-21,40){$(b_n,a_n)$}
\put(45,-5){$(r_{p},c_p)$}
\put(42,12){$(r_{p-1},c_{p-1})$}
\matrixput(0,0)(10,0){7}(0,10){4}{\circle*{.5}}
\matrixput(0,30)(10,0){4}(0,10){2}{\circle*{.5}}
\matrixput(0,40)(10,0){2}(0,10){2}{\circle*{.5}}
\linethickness{1pt}
\dottedline{0.1}(0,40)(0,30)(10,30)(10,20)(20,20)(20,10)(40,10)(40,0)(50,0)
\put(0,40){\makebox(0,0){$\bullet$}}
\put(10,30){\makebox(0,0){$\bullet$}}
\put(10,20){\makebox(0,0){$\bullet$}}
\put(20,20){\makebox(0,0){$\bullet$}}
\put(40,10){\makebox(0,0){$\bullet$}}
\put(50,0){\makebox(0,0){$\bullet$}}
\linethickness{0.1pt}
\multiput(0,0)(1.5,0){40}{\line(1,0){1}}
\multiput(0,10)(1.5,0){40}{\line(1,0){1}}
\multiput(0,20)(1.5,0){40}{\line(1,0){1}}
\multiput(0,30)(1.5,0){40}{\line(1,0){1}}
\multiput(0,40)(1.5,0){20}{\line(1,0){1}}
\multiput(0,50)(1.5,0){7}{\line(1,0){1}}

\multiput(0,0)(0,1.5){33}{\line(0,1){1}}
\multiput(10,0)(0,1.5){33}{\line(0,1){1}}
\multiput(20,0)(0,1.5){27}{\line(0,1){1}}
\multiput(30,0)(0,1.5){27}{\line(0,1){1}}
\multiput(40,0)(0,1.5){20}{\line(0,1){1}}
\multiput(50,0)(0,1.5){20}{\line(0,1){1}}
\multiput(60,0)(0,1.5){20}{\line(0,1){1}}
\put(76,40){$a_n$ bumps $c_1$ from the first row of $P^{(n-1)}$,}
\put(10.5,32){$(r_1,c_1)$}
\put(63,32) {$\longrightarrow$}
\put(76,32){$a_n$ is placed as an entry in the first row of $P^{(n)}$}
\put(76,24){and, $r_p$ is an entry in the first row of $Q^{(n)}$.}
\put(-5,-9.5){$\leftarrow$}
\put(-18,-8){$r_{p}$}
\multiput(0,-8)(1.15,0){44}{\line(1,0){1}}
\put(-3.1,-14){$\downarrow$}
\put(-5,-18){$a_n$}
\multiput(-2,-10)(0,1.15){42}{\line(0,1){1}}
\end{picture}

\begin{picture}(60,100)(0,0)
\put(-15,60){{Step II:}}
\put(-25,45){Case I :}
\matrixput(0,0)(10,0){7}(0,10){4}{\circle*{.5}}
\matrixput(0,30)(10,0){4}(0,10){2}{\circle*{.5}}
\matrixput(0,40)(10,0){2}(0,10){2}{\circle*{.5}}
\linethickness{1pt}
\dottedline{0.1}(10,40)(10,30)(20,30)(20,20)(40,20)(40,10)(50,10)
\put(-10,38){$(b_n,c_1)$}
\put(47,5){$(r_{p-1},c_p)$}
\put(10,40){\makebox(0,0){$\bullet$}}
\put(10,30){\makebox(0,0){$\bullet$}}
\put(20,20){\makebox(0,0){$\bullet$}}
\put(40,20){\makebox(0,0){$\bullet$}}
\put(50,10){\makebox(0,0){$\bullet$}}
\linethickness{0.1pt}
\multiput(0,0)(1.5,0){40}{\line(1,0){1}}
\multiput(0,10)(1.5,0){40}{\line(1,0){1}}
\multiput(0,20)(1.5,0){40}{\line(1,0){1}}
\multiput(0,30)(1.5,0){40}{\line(1,0){1}}
\multiput(0,40)(1.5,0){20}{\line(1,0){1}}
\multiput(0,50)(1.5,0){7}{\line(1,0){1}}

\multiput(0,0)(0,1.5){34}{\line(0,1){1}}
\multiput(10,0)(0,1.5){34}{\line(0,1){1}}
\multiput(20,0)(0,1.25){32}{\line(0,1){1}}
\multiput(30,0)(0,1.25){32}{\line(0,1){1}}
\multiput(40,0)(0,1.5){20}{\line(0,1){1}}
\multiput(50,0)(0,1.5){20}{\line(0,1){1}}
\multiput(60,0)(0,1.5){20}{\line(0,1){1}}
\put(76,50){$c_1$ (which was bumped by $a_n$ from the first row} 
\put(76,42){of $P^{(n-1)}$ in step I) is placed in the 2nd row of}
\put(76,34){the left hand tableau. And as $c_1$ is the least}
\put(76,26){column number of some block having at least two}
\put(76,18){elements, so it will bump an entry from the second row.}
\put(60,32) {$\longrightarrow$}
\put(-5,-9.5){$\leftarrow$}
\put(-18,-8){$r_{p-1}$}
\multiput(0,-8)(1.15,0){44}{\line(1,0){1}}
\put(6.68,-14){$\downarrow$}
\put(5,-20){$c_1$}
\multiput(8,-10)(0,1.15){42}{\line(0,1){1}}

\end{picture}

\begin{picture}(90,90)(0,0)
\put(-25,45){Case II :}
\put(-10,36.5){$\leftarrow$}
\put(-17,35){$b_n$}
\multiput(-5,38)(1.15,0){20}{\line(1,0){1}}
\put(20.59,-4){$\downarrow$}
\put(20,-10){$\star$}
\multiput(22,0)(0,1.15){34}{\line(0,1){1}}
\matrixput(0,0)(10,0){7}(0,10){4}{\circle*{.5}}
\matrixput(0,30)(10,0){4}(0,10){2}{\circle*{.5}}
\matrixput(0,40)(10,0){2}(0,10){2}{\circle*{.5}}
\linethickness{1pt}
\put(20,40){\makebox(0,0){$\bullet$}}

\put(18,44){$( b_n,\star)$}
\linethickness{0.1pt}
\multiput(0,0)(1.5,0){40}{\line(1,0){1}}
\multiput(0,10)(1.5,0){40}{\line(1,0){1}}
\multiput(0,20)(1.5,0){40}{\line(1,0){1}}
\multiput(0,30)(1.5,0){40}{\line(1,0){1}}
\multiput(0,40)(1.5,0){20}{\line(1,0){1}}
\multiput(0,50)(1.5,0){7}{\line(1,0){1}}

\multiput(0,0)(0,1.5){34}{\line(0,1){1}}
\multiput(10,0)(0,1.25){40}{\line(0,1){1}}
\multiput(20,0)(0,1.25){32}{\line(0,1){1}}
\multiput(30,0)(0,1.25){32}{\line(0,1){1}}
\multiput(40,0)(0,1.5){20}{\line(0,1){1}}
\multiput(50,0)(0,1.5){20}{\line(0,1){1}}
\multiput(60,0)(0,1.5){20}{\line(0,1){1}}
\put(60,32) {$\longrightarrow$}
\put(70,40){As $\{(b_n,\star)\}$ is a single block, so $b_n$ will take place}
\put(70,32){at the left end of some row in the right hand tableau}
\put(70,24){and $\star$ will take place at the right end of the}
\put(70,16){corresponding row in the left hand tableau.}
\end{picture}
\caption{An illustration of Lemma \ref{l.D}}
\label{figure.7}
\end{figure}
\begin{eg}\label{eg.lemmaD}
\newcommand{\A}{13}
\newcommand{\B}{12}
\newcommand{\M}{11}
\newcommand{\J}{10}
Let $d=7, N=13, v=(1,2,4,5,7,8,9)$. Let $$F=\{(10,4),(10,9),(11,1),(11,2),(11,7),(11,8),(12,2),(13,8)\}.$$
Then
\begin{equation}\label{equ.1}
BRSK(\iota(F))=(P^{(n-1)},Q^{(n-1)})= \left( \hspace{.2cm}\young(149,27,2,8,8)\hspace{.2cm} ,  \hspace{.2cm} \young(\J\M\A,\J\B,\M,\M,\M) \right). 
\end{equation}
\\
Figure~\ref{figure.8} below gives the block decomposition of $F$.\\
\begin{figure}[h]
\begin{picture}(80,70)(0,0)
\matrixput(0,0)(10,0){7}(0,10){4}{\circle*{.5}}
\matrixput(0,30)(10,0){4}(0,10){2}{\circle*{.5}}
\matrixput(0,40)(10,0){2}(0,10){2}{\circle*{.5}}
\linethickness{1pt}
\dottedline{0.1}(0,20)(10,20)(10,0)(50,0)
\put(0,20){\makebox(0,0){$\bullet$}}
\put(10,20){\makebox(0,0){$\bullet$}}
\put(10,10){\makebox(0,0){$\bullet$}}
\put(50,0){\makebox(0,0){$\bullet$}}
\put(20,30){\makebox(0,0){$\bullet$}}
\put(40,20){\makebox(0,0){$\bullet$}}
\put(50,20){\makebox(0,0){$\bullet$}}
\put(60,30){\makebox(0,0){$\bullet$}}
\linethickness{.9pt}
\dottedline{0.1}(20,30)(20,20)(50,20)
\linethickness{0.1pt}
\multiput(0,0)(1,0){60}{\line(1,0){1}}
\multiput(0,10)(1,0){60}{\line(1,0){1}}
\multiput(0,20)(1,0){60}{\line(1,0){1}}
\multiput(0,30)(1,0){60}{\line(1,0){1}}
\multiput(0,40)(1,0){30}{\line(1,0){1}}
\multiput(0,50)(1,0){10}{\line(1,0){1}}
\linethickness{0.27pt}
\multiput(0,0)(0,1){50}{\line(0,1){1}}
\multiput(10,0)(0,1){50}{\line(0,1){1}}
\multiput(20,0)(0,1){40}{\line(0,1){1}}
\multiput(30,0)(0,1){40}{\line(0,1){1}}
\multiput(40,0)(0,1){30}{\line(0,1){1}}
\multiput(50,0)(0,1){30}{\line(0,1){1}}
\multiput(60,0)(0,1){30}{\line(0,1){1}}
\multiputlist(-10,0)(0,10){13,12,11,10,6,3}
\put(0,55){1}
\put(10,55){2}
\put(20,45){4}
\put(30,45){5}
\put(40,35){7}
\put(50,35){8}
\put(60,35){9}
\end{picture}
\caption{This is the block decomposition of the monomial $F$. Here $F_1=\{(11,1),(11,2),(12,2),(13,8)\}$, $F_2=\{(10,4),(11,7),(11,8)\}$, and $F_3=\{(10,9)\}$, where $F_j$ is the set of all elements of $F$ which are $j$ deep but not $(j+1)$ deep. }
\label{figure.8}
\end{figure}
Now let $(b_n,a_n)=(6,2)$. Hence 
$$U=\{(6,2),(10,4), (10,9),(11,1),(11,2),(11,7),(11,8),(12,2),(13,8)\}.$$
Then
\begin{equation}\label{equ.2}
BRSK(\iota(U))=(P^{(n)},Q^{(n)})= \left( \hspace{.2cm}\young(129,247,2,8,8)\hspace{.2cm} ,  \hspace{.2cm} \young(\J\M\A,6\J\B,\M,\M,\M) \right).
\end{equation}
Figure~\ref{figure.9} below gives the block decomposition of $U$.\\
\begin{figure}[h]
\begin{picture}(50,60)(0,0)
\matrixput(0,0)(10,0){7}(0,10){4}{\circle*{.5}}
\matrixput(0,30)(10,0){4}(0,10){2}{\circle*{.5}}
\matrixput(0,40)(10,0){2}(0,10){2}{\circle*{.5}}
\linethickness{1pt}
\dottedline{0.1}(0,20)(10,20)(10,10)(10,0)(50,0)
\put(0,20){\makebox(0,0){$\bullet$}}
\put(10,20){\makebox(0,0){$\bullet$}}
\put(10,10){\makebox(0,0){$\bullet$}}
\put(10,40){\makebox(0,0){$\bullet$}}
\put(50,0){\makebox(0,0){$\bullet$}}
\put(20,30){\makebox(0,0){$\bullet$}}
\put(40,20){\makebox(0,0){$\bullet$}}
\put(50,20){\makebox(0,0){$\bullet$}}
\put(60,30){\makebox(0,0){$\bullet$}}
\linethickness{1pt}
\dottedline{0.1}(10,40)(10,30)(20,30)(20,20)(50,20)
\linethickness{0.1pt}
\multiput(0,0)(1,0){60}{\line(1,0){1}}
\multiput(0,10)(1,0){60}{\line(1,0){1}}
\multiput(0,20)(1,0){60}{\line(1,0){1}}
\multiput(0,30)(1,0){60}{\line(1,0){1}}
\multiput(0,40)(1,0){30}{\line(1,0){1}}
\multiput(0,50)(1,0){10}{\line(1,0){1}}
\linethickness{0.4pt}
\multiput(0,0)(0,1.25){40}{\line(0,1){1}}
\multiput(10,0)(0,1.25){40}{\line(0,1){1}}
\multiput(20,0)(0,1.25){32}{\line(0,1){1}}
\multiput(30,0)(0,1.25){32}{\line(0,1){1}}
\multiput(40,0)(0,1.5){20}{\line(0,1){1}}
\multiput(50,0)(0,1.5){20}{\line(0,1){1}}
\multiput(60,0)(0,1.5){20}{\line(0,1){1}}
\multiputlist(-10,0)(0,10){13,12,11,10,6,3}
\put(0,55){1}
\put(10,55){2}
\put(20,45){4}
\put(30,45){5}
\put(40,35){7}
\put(50,35){8}
\put(60,35){9}
\end{picture}
\caption{This is the block decomposition of the monomial $U$. Here $U_1=\{(11,1),(11,2),(12,2),(13,8)\}$, $U_2=\{(6,2),(10,4),(11,7),(11,8)\}$, and $U_3=\{(10,9)\}$, where $U_j$ is the set of all elements of $U$  which are $j$ deep but not $(j+1)$ deep.}
\label{figure.9}
\end{figure}
\vspace{2 cm}

So $U^{(1)}=\{(11,2),(11,2),(12,8),(6,4),(10,7),(11,8)\}$. Figure~\ref{figure.10} below gives the block decomposition of $U^{(1)}$.\\
\begin{figure}[h]
\begin{picture}(60,50)(0,0)
\matrixput(0,0)(10,0){7}(0,10){4}{\circle*{.5}}
\matrixput(0,30)(10,0){4}(0,10){2}{\circle*{.5}}
\matrixput(0,40)(10,0){2}(0,10){2}{\circle*{.5}}
\linethickness{1pt}
\dottedline{0.1}(10,20)(50,20)(50,10)
\put(10,20){\makebox(0,0){$\bullet$}}
\put(50,20){\makebox(0,0){$\bullet$}}
\put(20,40){\makebox(0,0){$\bullet$}}
\put(40,30){\makebox(0,0){$\bullet$}}
\put(50,10){\makebox(0,0){$\bullet$}}
\linethickness{0.1pt}
\multiput(0,0)(1,0){60}{\line(1,0){1}}
\multiput(0,10)(1,0){60}{\line(1,0){1}}
\multiput(0,20)(1,0){60}{\line(1,0){1}}
\multiput(0,30)(1,0){60}{\line(1,0){1}}
\multiput(0,40)(1,0){30}{\line(1,0){1}}
\multiput(0,50)(1,0){10}{\line(1,0){1}}
\linethickness{0.3pt}
\multiput(0,0)(0,.75){67}{\line(0,1){1}}
\multiput(10,0)(0,1){50}{\line(0,1){1}}
\multiput(20,0)(0,1){40}{\line(0,1){1}}
\multiput(30,0)(0,1){40}{\line(0,1){1}}
\multiput(40,0)(0,1){30}{\line(0,1){1}}
\multiput(50,0)(0,1){30}{\line(0,1){1}}
\multiput(60,0)(0,1){30}{\line(0,1){1}}
\multiputlist(-10,0)(0,10){13,12,11,10,6,3}
\put(0,55){1}
\put(10,55){2}
\put(20,45){4}
\put(30,45){5}
\put(40,35){7}
\put(50,35){8}
\put(60,35){9}
\end{picture}
\caption{This is the block decomposition of the monomial $U^{(1)}$. Here $U_1^{(1)}=\{(11,2),(11,2),(11,8),(12,8)\}$. $U_{2}^{(1)}$ is a union of two blocks namely, $\{(6,4)\}$ and $\{(10,7)\}$. Here $U_{j}^{(1)}$ is the set of all elements of $U^{(1)}$ which are $j$ deep but not $(j+1)$ deep.}
\label{figure.10}
\end{figure}
\vspace{.002 cm}

Here $b_n=6, a_n=2$. From Figure~\ref{figure.8}, we can see that $\{(10,4),(11,7),(11,8)\}$ is the topmost block of $F$ of depth $k=2$, and from Figure~\ref{figure.9}, we can see that after adding $(6,2)$ to $F$, $\{(6,2),(10,4),(11,7),(11,8)\}$ becomes the topmost block of $U$ of depth $k=2$. From Equation \ref{equ.2}, we can see that $a_n=2$ is an entry in the first row of $P^{(n)}$ and $r_p=11$ is an entry of the first row of $Q^{(n)}$. $a_n$ bumps $c_1=4$ from the first row of $P^{(n-1)}$. $c_1$ is placed in the second row of $P^{(n-1)}$. Also from Figure~\ref{figure.10}, we can see that since $\{(b_n, \star)\}=\{(6,4)\}$ is a single block, so $\star=4$ is placed in a new box in the second row of $P^{(n-1)}$, and $b_n$ is placed in a new box in second row of $Q^{(n-1)}$. 
\end{eg}
\begin{Proof}
 By induction hypothesis, the first row of $P^{(n-1)}$ contains the smallest column number of each block of $F$. In particular, it contains the entry $c_1$. And the first row of $Q^{(n-1)}$ contains the entry $r_p$.
In the process $P^{(n-1)}\la{b_n} a_n$ of bounded insertion, firstly all entries of $P^{(n-1)}$ which are $\geq b_n$ are removed. Since $c_1 < b_n$, $c_1$ is not removed.\par
\noindent\textbf{Claim:} $a_n$ bumps $c_1$ from the first row of $P^{(n-1)}$.\\
\noindent\textbf{Proof of claim:} Clearly $a_n \leq c_1$. It suffices to show that all the entries in the first row of $P^{(n-1)}$ which are $<c_1$ are also strictly less than $a_n$.\par 
Suppose not. Say, $p_{1j_0}$ is an entry in the first row of $P^{(n-1)}$ such that $a_n \leq p_{1j_0} < c_1$. By Lemma \ref{l.A}, it follows that $p_{1j_0}$ is the column number of the topmost element of some block $\mathfrak{C}$ of $F$ of depth $s (<k)$. Say, the first (or the topmost) element of the block $\mathfrak{C}$ is $(R,p_{1j_0})$. Now, since $(b_n,a_n)$ is an element of depth $k$ in $U$ and $s<k$, there exists an element $(e_s,f_s)$ in $U$ of depth $s$ such that $(e_s,f_s) > (b_n,a_n)$. Clearly then  $e_s >b_n$ and $f_s<a_n \leq p_{1j_0}$. Also $e_s \leq R$, because if $e_s>R$, then $(e_s,f_s) > (R,p_{1j_0})$, a contradiction to the depth of $(R,p_{1j_0})$ in $U$.\par 
Hence we have $b_n < e_s \leq R$ and $f_s < p_{1j_0}$. Now $(e_s,f_s)$ and $(R, p_{1j_0})$ are two elements in $U$ of depth $s$. Also $e_s > b_n > c_1 > p_{1j_0}$. Which implies $e_s > p_{1j_0}$, that is, $(e_s,f_s)$ belongs to the same block $\mathfrak{C}$ as $(R,p_{1j_0})$. Since $e_s \leq R$ and $f_s < p_{1j_0}$, we get a contradiction to the fact that $(R,p_{1j_0})$ is the topmost element of the block $\mathfrak{C}$.\par 
This proves the above claim.\par 
Now, since $a_n$ bumps $c_1$ from the first row of $P^{(n-1)}$, it follows that $a_n$ is placed at the position of $c_1$ in the first row of $P^{(n)}$. Also $c_1$ gets inserted in the 2nd row of $P^{(n-1)}$. Again, it is clear that $\exists$ a block of $U^{(1)}$ whose smallest column number is $c_1$ (the proof of this fact follows from the proof of Lemma~\ref{l.C} above). If $c_1$ does not bump anything from the second row of $P^{(n-1)}$, that means that $\{(b_n,c_1)\}$ is itself a block of $U^{(1)}$. In this case, $b_n$ is placed in a new box in the 2nd row of $Q^{(n-1)}$ and the process $(P^{(n-1)},Q^{(n-1)})\la{b_n}a_n$ terminates. If $c_1$ bumps something (say, $p_{2j_0}$) from the 2nd row of $P^{(n-1)}$, then $c_1 \leq p_{2j_0}$ and the 2nd row of $Q^{(n-1)}$ remains unchanged.\par
Now the remaining statements of the lemma follow from Lemma \ref{l.C}.
\end{Proof}
\section{Proof of Theorem \ref{t.main}}\label{s.theproof}
\begin{Proof}
Let degree of $U$ be $n$. We will prove the theorem by induction on $n$. Suppose the theorem is true for all monomials in $\widetilde{\mathfrak{N}^v}$ of degree $\leq n-1$. Arrange $ \iota(U)$ in lexicographic order, say, $\iota(U)=\{(a_1,b_1),(a_2,b_2),\ldots,(a_n,b_n)\}$. Let $F=\{(b_1,a_1),(b_2,a_2),\ldots,(b_{n-1},a_{n-1})\}$. That is $\iota(F)=\{(a_1,b_1),(a_2,b_2),\ldots,(a_{n-1},b_{n-1})\}$. Then $b_n \leq b_{n-1}$ and if $b_n = b_{n-1}$, we have $a_{n-1} \geq a_n$. The following cases arise:\\
(i) $b_n< b_{n-1} $ and $a_{n-1} > a_n$.\\
(ii) $b_n< b_{n-1} $ and $a_{n-1} = a_n$.\\
(iii) $b_n< b_{n-1} $ and $a_{n-1} < a_n$.\\
(iv) $b_n= b_{n-1} $ and $a_{n-1} = a_n$.\\
(v) $b_n= b_{n-1} $ and $a_{n-1} > a_n$.\\
$(b_n,a_n)$ enters into $F$ to make it $U$.\\
\noindent\textbf{Case (i):} $b_n< b_{n-1} $ and $a_{n-1} > a_n$:\\
Since $b_n< b_{n-1} $ , so $(b_n,a_n)$ cannot change the depth of any element of $F$. Again since $a_{n-1}>a_n$ and the elements in $\{(a_1,b_1),\ldots,(a_n,b_n)\}$ are in lexicographic order, it follows that $(b_n,a_n)$ cannot make a new block of higher depth (that is, if the maximum possible depth of any element of $F$ is $t$, then it is not possible that the maximum possible depth of any element in $U$ becomes strictly bigger than $t$). So the only thing that can happen is that the element $(b_n,a_n)$ gets added to $F$ as the topmost element of some block. Hence after $(b_n,a_n)$ enters into $F$, only two things can happen :
Either $\{(b_n,a_n\})$ becomes the topmost block of $U$ of depth $k$ containing a single element or $\{(b_n,a_n),(r_1,c_1),\ldots,(r_p,c_p)\}$ becomes the topmost block $U$ of depth $k$, where $\{(r_1,c_1),\ldots,(r_p,c_p)\}$ is the topmost block of $F$ of depth $k$.\par
\noindent\textbf{Subcase (a):} Suppose the singleton set $\{(b_n,a_n)\}$ becomes the topmost block of $U$ of depth $k$ :\\
$\{(b_n,a_n)\}$ is the topmost block of $U$ of depth $k$, the next block of $U$ of depth $k$ from the top being $\{(r_1,c_1),(r_2,c_2),\ldots,(r_p,c_p)\}$. Then $b_n<c_1, \ b_n<r_1 \leq r_2 \leq\ldots\leq r_p$ and $a_n < c_1 \leq c_2 \leq\ldots\leq c_p$. Let $BRSK(\iota(F))=(P^{(n-1)},Q^{(n-1)})$. Let $p_{ij}$ denote the entry in the $i$-th row and $j$-th column of $P^{(n-1)}$. Similarly $q_{ij}$. Clearly, $BRSK(\iota(U))= (P^{n-1},Q^{n-1})\la{b_n}a_n$. By induction hypothesis, it follows that the first row of $P^{(n-1)}$ contains the smallest column numbers of each block of $F$. In particular, it contains the entry $c_1$. Then by Lemma \ref{l.B}, we know that the entries in the first row of $P^{(n-1)}$ which are strictly less than $b_n$ are all strictly less than $a_n$. Hence we are done.\par
\noindent\textbf{Subcase (b):} $\{(b_n,a_n),(r_1,c_1),\ldots,(r_p,c_p)\}$ becomes the topmost block of $U$ of depth $k$.\\
Then we have $b_n > c_1, \  b_n<r_1 \leq r_2 \leq\ldots\leq r_p$ and $a_n \leq c_1 \leq c_2 \leq\ldots\leq c_p$.\\
Let $(P^{(n)},Q^{(n)})$ = $BRSK(\iota(U))=(P^{(n-1)},Q^{(n-1)}) \la{b_n}a_n$. By Lemma \ref{l.D} we are done.\\
\noindent\textbf{Case (ii):} $b_n < b_{n-1}$ and $a_{n-1}=a_n$:\\
Since $(b_1,a_1),\ldots,(b_n,a_n)$ are in lexicographic order, it follows that $(b_{n-1},a_{n-1})$ is the topmost element of the topmost block of $F$ of some depth, say $k$. In this case, since $(b_1,a_1),\ldots,(b_n,a_n)$ are in lexicographic order, $b_n<b_{n-1}$ and $a_{n-1}=a_n$, it follows that $(b_n,a_n)$ can only get added to the block of $F$ whose topmost element is $(b_{n-1},a_{n-1})$. Also since $b_n < b_{n-1}$ and $a_{n-1}=a_n$, $(b_{n},a_{n})$ will be the topmost element of that block. Now if $\{(b_{n-1},a_{n-1})\}$ had been a single block of $F$ of depth $k$, then $\{(b_n,a_n),(b_{n-1},a_{n-1})\}$ will be the first block of $U$ of depth $k$. Then the proof of the theorem is similar to the proof of Subcase (b), Case (i). And if $\{(b_{n-1},a_{n-1}),(r_1,c_1),\ldots,(r_p,c_p)\}$ is the topmost block of $F$ of depth $k$, then $\{(b_n,a_n),(b_{n-1},a_{n-1}),(r_1,c_1),\ldots,(r_p,c_p)\}$ will be the topmost block of $U$ of depth $k$. Then again the proof of the theorem is similar to the proof of subcase (b), case(i).\\
\noindent\textbf{Case (iii):} $b_n <b_{n-1}$ and $a_{n-1} < a_n$:\\
Since $b_n<b_i \  \forall i \in \{1,\ldots,n-1\}$, $(b_n,a_n)$ cannot change the depth of any element  of $F$. So only two things can happen, which are given in the Subcases (a) and (b) below:\\
\noindent\textbf{Subcase (a):} $(b_n,a_n)$ gets added to a block of $F$ for some depth $k$ as the topmost element.\\
\noindent\textbf{Proof:} Same as in Case (i).\\
\noindent\textbf{Subcase (b):} The singleton set $\{(b_n,a_n)\}$ becomes the topmost block of depth $(k+1)$ in $U$, where $k$ is the maximum possible depth of any element of $F$.\\
\noindent\textbf{Proof:} In this subcase, we will have to show that $a_n$ is strictly greater than all entries in the first row of $P^{(n-1)}$ which are $< b_n$.\\
Suppose not.\\
Then $\exists$ an entry (call it $p_{1j_0}$) in the first row of $P^{(n-1)}$ which is such that
\begin{equation}\label{equ.3}
a_n \leq p_{1j_0} <b_n.
\end{equation}
Clearly then, $p_{1j_0}$ is the smallest column number of some block $\mathfrak{B}$ of $F$, say of depth $s$, where $s\leq k$.\\
Let $(R,p_{1j_0})$ denote the topmost element of the block $\mathfrak{B}$.\\
\noindent\textbf{Claim 1 :} The block $\mathfrak{B}$ cannot be the topmost block of depth $s$ in $F$.\\
\textbf{Proof of Claim 1 :} Suppose not. Say $\mathfrak{B}$ is the topmost block of depth $s$ in $F$. Since $(b_n,a_n)$ has depth $k+1$ in $U$ and $s < k+1$, there exists an element $(\beta,\alpha)$ such that $(\beta,\alpha)$ has depth $s$ in $F$ and $(\beta,\alpha)>(b_n,a_n)$ is a $v$-chain. Then $b_n < \beta$ and $a_n>\alpha$. The element $(\beta,\alpha)$ of $F$ lies in some block of $F$ of depth $s$. Since $\mathfrak{B}$ is the topmost block of depth $s$ in $F$ and $p_{1j_0}$ is the smallest column number of $\mathfrak{B}$, we must have $p_{1j_0} \leq \alpha$. But then we get $p_{1j_0} \leq \alpha <a_n$, which contradicts equation (\ref{equ.3}). Hence Claim 1 is proved.\\
\noindent\textbf{Claim 2 :} $\exists$ an element $(e_s,f_s)$ in the topmost block of depth $s$ in $F$ such that $(e_s,f_s)>(b_n,a_n)$ is a $v$-chain.\\
\textbf{Proof of Claim 2 :} Since $(b_n,a_n)$ is of depth $k+1$ in $U$ and $s<k+1$, there exists an element $(e_s,f_s)$ of depth $s$ in $F$ such that $(e_s,f_s)>(b_n,a_n)$ is a $v$-chain. It suffices to show that $(e_s,f_s)$ belongs to the topmost block of depth $s$ in $F$.\par 
Suppose not. Let $\mathfrak{B}$ denote the block of $F$ of depth $s$ where $(e_s,f_s)$ lies.\par 
Let $\mathfrak{C}$ denote the topmost block of $F$ of depth $s$. Let $(b,a)$ denote the bottom-most element of $\mathfrak{C}$. Then since $\mathfrak{B}\neq \mathfrak{C}$, we have $b \leq f_s$. Now since $(e_s,f_s) > (b_n,a_n)$, we have $f_s < a_n$. So we get $b \leq f_s < a_n$, which implies $(b,a_n) \notin\widetilde{\mathfrak{N}^v}$.    But $b_n < b$ (since $b_n < b_{i}$ for all $i\in\{1,\ldots,n-1\}$). So, we get $b_n < a_n$, a contradiction to the fact that $(b_n,a_n) \in\widetilde{\mathfrak{N}^v}$. Hence Claim 2 is proved.\\
\noindent\textbf{Assuming Claim 1 and Claim 2:} It follows from Claim 2 that $\exists\ (e_s,f_s)$ in the topmost block of depth $s$ in $F$ such that $(e_s,f_s) > (b_n, a_n)$ is a $v$-chain. $\Rightarrow b_n < e_s$. Now, Equation \ref{equ.3} says that $a_n \leq p_{1j_0} < b_n$. This together with $b_n < e_s$ imply that $p_{1j_0}< e_s$. But  this implies the two elements $(e_s,f_s)$ and $(R,p_{1j_0})$ of depth $s$ lie in the same block, which is a contradiction to claim 1.\\
\noindent \textbf{Case (iv):} $b_n = b_{n-1}$ and $a_n = a_{n-1}$: \\
In this case $(b_n,a_n)$ will be added in the block of $(b_{n-1},a_{n-1})$ and it will be the topmost element of that block. The rest follows similarly as in Subcase (b) of Case (i).\\
\noindent\textbf{Case (v):} $ b_n = b_{n-1}$ and $a_{n-1} > a_n$: \\
Since $ b_n = b_{n-1}$ and $a_{n-1} > a_n$, so $(b_n,a_n)$ cannot make a new block of higher depth. Again since $b_n \leq b_i \  \forall \ i \in \{1,\ldots,n-1\}$, it follows that $(b_n,a_n)$ cannot change the depth of any element of $F$. So the only thing that can happen is that the element $(b_n,a_n)$ gets added to $F$ as the topmost element of some block (of some depth $k$ say). In this case also, the proof is similar to the proof of Case(i).
\end{Proof}
\section{An application of the main theorem}\label{s.theproblem}
\subsection{Some necessary definitions and notation}\label{ss.notation}
\noindent The following definitions and notation are written in the same way as given in the papers \cite{gr} and \cite{kr}.\par 
A positive integer $d$ will be kept fixed throughout this section. For $j \in [2d]$, set $j^* : = 2d+1-j$. Let $I(d)$ denote the set of all $d$-element subsets $v$ of $[2d]$ with the property that exactly one of $j$, $j^*$ belongs to $v$ for every $j \in [d]$. Clearly $I(d)\subseteq I(d,2d)$. In particular, we have the partial order $\leq$ on $I(d)$ induced from $I(d,2d)$. We denote by $\epsilon$ the element $(1,\ldots,d)$ of $I(d)$. The $\epsilon$-\textbf{degree} of an element $x$ of $I(d)$ is the cardinality of $x \setminus [d]$ or equivalently that of $[d]\setminus x$. More generally, given any $v\in I(d)$, the $v$-\textbf{degree} of an element $x$ of $I(d)$ is the cardinality of $x\setminus v$ or equivalently that of $v\setminus x$. 
\begin{eg}\label{exmple.1}
Let $d=5$. So $2d=10$, and $\epsilon =(1,2,3,4,5)$. Let $x=(1,2,4,6,8)$. Hence in this case, the $\epsilon$-degree of $x$ is $2$. Again let $v=(1,2,3,4,6)$. So according to the definition, $v$-degree of $x$ is $1$. 
\end{eg}
An ordered pair $\mathfrak{w}=(x,y)$ of elements of $I(d)$ is called an \textbf{admissible pair} if $x\geq y$ and the $\epsilon$-degrees of $x$ and $y$ are equal. Sometimes $x$ and $y$ are referred as the top and bottom of $\mathfrak{w}$ and written as $top(\mathfrak{w})$ for $x$ and $bot(\mathfrak{w})$ for $y$. Given an admissible pair $\mathfrak{w}=(x,y)$, we define the $v$-\textbf{degree} of $\mathfrak{w}$ by $v$-degree$(\mathfrak{w}):=\frac{1}{2}(|x\setminus v|+|y\setminus v|)$.\par 
Given any two admissible pairs $\mathfrak{w}=(x,y)$ and $\mathfrak{w'}=(x',y')$, we say that $\mathfrak{w} \geq \mathfrak{w'}$ if $y\geq x'$, that is , if $x\geq y \geq x' \geq y'$. An ordered sequence $(\mathfrak{w}_1 \ldots, \mathfrak{w}_t)$ of admissible pairs is called a \textbf{standard sequence} if $\mathfrak{w}_i\geq \mathfrak{w}_{i+1}$ for $1\leq i<t$. Sometimes $\mathfrak{w}_1 \geq \ldots \geq \mathfrak{w}_t$ is written to denote a standard sequence $(\mathfrak{w}_1 \ldots, \mathfrak{w}_t)$. Given any $w\in I(d)$, we say that a standard sequence $\mathfrak{w}_1 \geq \ldots \geq \mathfrak{w}_t$ is $w$-\textbf{dominated} if $w\geq top(\mathfrak{w}_1)$. Given any $v\in I(d)$, we say that the standard sequence $\mathfrak{w}_1 \geq \ldots \geq \mathfrak{w}_t$ is $v$-\textbf{compatible} if for each $\mathfrak{w}_i$, either $v\geq top(\mathfrak{w}_i)$ or $bot(\mathfrak{w}_i)\geq v$, and $\mathfrak{w}_i\neq (v,v)$. Given $v$ and $w$ in $I(d)$, we denote by $SM_w^v$ the set of all $w$-dominated $v$-compatible standard sequences. For any positive integer $m$, let $SM_w^v(m)$ denote the set of all $w$-dominated $v$-compatible standard sequences of degree $m$, where the \textbf{degree of a standard sequence} $\mathfrak{w}_1 \geq \ldots \geq \mathfrak{w}_t$ is defined to be the sum of the $v$-degrees of $\mathfrak{w}_1,\ldots,\mathfrak{w}_t$. Let $SM^{v,v}$ denote the set of all $v$-compatible standard sequences that are anti-dominated by $v$: a standard sequence $\mathfrak{w}_1 \geq \ldots \geq \mathfrak{w}_t$ is called \textbf{anti-dominated} by $v$ if $bot(\mathfrak{w}_t)\geq v$. 
\begin{eg}\label{exmple.2}
Let $d=5$. So $2d=10$, and $\epsilon =(1,2,3,4,5)$. Let $x_1 =(3,4,6,9,10)$, and $y_1 =(2,4,6,8,10)$. Clearly, the $\epsilon$-degrees of $x_1$ and $y_1$ are equal (which is $3$), and $x_1 \geq y_1$. Hence $\mathfrak{w_1} = (x_1,y_1)$ is an admissible pair.\par 
Let $\mathfrak{w_2}=(x_2,y_2)$, where $x_2=(1,3,6,7,9)$, $y_2=(1,2,6,7,8)$. Let $\mathfrak{w_3}=(x_3,y_3)$, where $x_3=(1,2,5,7,8)$, $y_3=(1,2,4,6,8)$. Then $\mathfrak{w_2}$ and $\mathfrak{w_3}$ are both admissible pairs such that $\mathfrak{w_1} \geq \mathfrak{w_2} \geq \mathfrak{w_3}$. Hence $\mathfrak{w_1} \geq \mathfrak{w_2} \geq \mathfrak{w_3}$ is a standard sequence.\par 
Let $w=(4,5,8,9,10)$. Then clearly $w \geq top(\mathfrak{w_1})$. Hence $\mathfrak{w_1} \geq \mathfrak{w_2} \geq \mathfrak{w_3}$ is dominated by $w$. Again for $v=(1,2,3,4,6) \in I(d)$, $bot(\mathfrak{w_i}) \geq v \ \forall \ i\in \{1,2,3\}$, and $\mathfrak{w_i} \neq (v,v)$ for all $i\in\{1,2,3\}$. So the above standard sequence is $v$-compatible. Hence it belongs to $SM_w^v$.\par 
Also $bot(\mathfrak{w_3}) \geq v$, so $\mathfrak{w_1} \geq \mathfrak{w_2} \geq \mathfrak{w_3}$ is anti-dominated by $v$. So the standard sequence belongs to $SM^{v,v}$.
\end{eg}
Fix a vector space $V$ of dimension $2d$ over an algebraically closed field of arbitrary characteristic. Fix a \textbf{non degenerate skew-symmetric bilinear form} $\langle\ ,\ \rangle$ on $V$. Fix a basis $e_1,\ldots,e_{2d}$ of $V$ such that \par 
$\langle e_i,e_j\rangle$=$\left\{\begin{array}{lr}
1 \hspace{.5cm} \mbox{if} \hspace{.2cm} i=j^* \hspace{.2cm} \mbox{and} \hspace{.2cm} i<j \\
-1 \hspace{.5cm} \mbox{if} \hspace{.2cm} i=j^*\hspace{.2cm} \mbox{and} \hspace{.2cm}i>j \\
0 \hspace{2cm} \mbox{otherwise}.
\end{array}\right.$\par 
A linear subspace $W$ of $V$ is said to be \textbf{isotropic} if the form $\langle\ ,\ \rangle$ \textbf{vanishes identically} on it. Denote by $\mathfrak{G}_d(V)$ the Grassmannian of all $d$-dimensional subspaces of $V$ and by $\mathfrak{M}_d(V)$ the set of all maximal isotropic subspaces of $V$. Then $\mathfrak{M}_d(V)$ is a closed subvariety of $\mathfrak{G}_d(V)$ and is called the \textbf{symplectic Grassmannian}.\par 
The group $Sp(V)$ of linear automorphisms of $V$ preserving $\langle\ ,\ \rangle$ acts transitively on $\mathfrak{M}_d(V)$-this follows from Witt's theorem that an isometry between subspaces can be lifted to one of the whole vector space. The elements of $Sp(V)$ that are diagonal with respect to the basis $\{e_1, \ldots ,e_{2d}\}$ form a maximal torus $T$ of $Sp(V)$. Similarly the elements of $Sp(V)$ that are upper triangular with respect to $\{e_1, \ldots ,e_{2d}\}$ form a Borel subgroup $B$ of $Sp(V)$.\par The $T$-fixed points of $\mathfrak{M}_d(V)$ are parametrized by $I(d)$: for $v=(v_1,\ldots,v_d)$ in $I(d)$, the corresponding $T$-fixed point, denoted by $e^v$, is the span of $e_{v_1},\ldots, e_{v_d}$. These points lie in different $B$-orbits, and the union of their $B$-orbits is all of $\mathfrak{M}_d(V)$. A \textbf{Schubert variety} in $\mathfrak{M}_d(V)$ is by definition the closure of such a $B$-orbits with the reduced scheme structure. Schubert varieties are thus indexed by the $T$-fixed points and so in turn by $I(d)$. Given $w$ in $I(d)$, we denote by $X_w$ the closure of the $B$-orbit of the $T$-fixed point $e^w$.\par 
Fix an elements $v$ of $I(d)$. Define $\mathfrak{R}^v$ $ :=\{(r,c) \in [2d] \setminus v \times v : r \leq c^\ast\}$ and let $\widetilde{\mathfrak{N}^v}:=\{(r,c) \in [2d] \setminus v \times v : r > c\}$. Let $S^v$ denote the set of all monomials in $\mathfrak{R}^v$ and $\widetilde{T^v}$ the set of all monomials in $\widetilde{\mathfrak{N}^v}$. Let $w$ be another element of $I(d)$ such that $v\leq w$. Let $S^v_w$ denote the set of all $w$-dominated monomials in $\mathfrak{R}^v$ (where $w$-domination of a monomial is defined as in Subsection \ref{ss.recap}), and (for any positive integer $m$) let $S^v_w(m)$ denote the set of such monomials of degree $m$.
\begin{eg}\label{exmple.3}
Let $d=5$, and $v=(1,2,4,6,8)$. Clearly $v$ is in $I(d)$. Let $\mathfrak{S}=\{(3,4)^2,(3,6),(5,4),(7,2),(9,1)\}$. So $\mathfrak{S}$ is a monomial in $\mathfrak{R}^v$. Now degree of $\mathfrak{S}$ is $6$, and $\beta_1 =(9,1) > \beta_2 =(7,2) > \beta_3=(5,4)$ is a $v$-chain in $\mathfrak{S}$. Again $s_{\beta_1}\cdots s_{\beta_3}v=(\{1,2,4,6,8\} \setminus \{1,2,4\}) \cup \{5,7,9\}=(5,6,7,8,9)$. Let $w=(5,7,8,9,10)$. The $v$-chain $\beta_1>\beta_2>\beta_3$ is dominated by $w$. Similarly, we can check that any other $v$-chain in $\mathfrak{S}$ is also dominated by $w$. Hence $\mathfrak{S}$ is dominated by $w$. So $\mathfrak{S}$ is an element of $S_w^v(6)$ .
\end{eg}
For $\alpha=(r,c)$ in $\mathfrak{R}^v$, set $\alpha ^\# :=(c^\star,r^\star)$. Elements of the form $(r,r^\star)$ of $\widetilde{\mathfrak{N}^v}$ are referred to as belonging to the ``diagonal", and the set of all diagonal elements of $\widetilde{\mathfrak{N}^v}$ is denoted by $\mathfrak{d}^v$. A monomial $\mathfrak{S}$ of $\widetilde{T^v}$ is \textbf{special} if\\ 
(1) $\mathfrak{S}=\mathfrak{S}^\#$ and\\
(2) the multiplicity of every diagonal element in $\mathfrak{S}$ is even.\\
Equivalently, $\mathfrak{S}$ is special if there exists $\mathfrak{T}$ in $\widetilde{T^v}$ with $\mathfrak{S}=\mathfrak{T} \cup \mathfrak{T}^\#$. The set of all special monomials is denoted by $\mathfrak{E}$.
\begin{eg}\label{exmple.4}
Let $d=5$, and $v=(1,2,4,6,8)$. Let $\mathfrak{T}=\{(7,3),(6,5)\}$. So $\mathfrak{T}^\#=\{(8,4),(6,5)\}$. Hence $\mathfrak{S}= \mathfrak{T} \cup \mathfrak{T}^\# = \{(6,5)^2,(7,3),(8,4)\}$ is a special monomial.
\end{eg}
\subsection{The result of Ghorpade and Raghavan}\label{ss.gr}
\begin{theorem}\label{t.gr}
Let $v$, $w$ be elements of $I(d)$ with $v\leq w$. Let $X_w$ be the Schubert variety corresponding to $w,e^v$ the $T$-fixed point in $X_w$ corresponding to $v$, and $R$ be the coordinate ring of the tangent cone to $X_w$ at the point $e^v$. Then the dimension as a vector space of the $m^{th}$ graded piece $R(m)$ of $R$ equals the cardinality of $S_w ^v(m)$.
\end{theorem}
\subsection{The application}\label{ss.ourproblem}
Before going to the application, we first need to state Proposition $4.1$ of \cite{gr}, which is given below as Proposition~\ref{prop.4.1}.
\begin{prop}\label{prop.4.1}
There is a bijection between $SM^{v,v}$ and $\widetilde{T^v}$ that respects domination and degree.
\end{prop}
The proof of Theorem \ref{t.gr} (as given in \cite{gr}) relies on a bijection between the two combinatorially defined sets $SM^v_w(m)$ and $S^v_w(m)$. And this bijection in turn, relies upon a bijection between $SM^{v,v}$ and $\widetilde{T^v}$, which is stated in Proposition~\ref{prop.4.1} above. The application of the main theorem here is to prove that the bijection between $SM^{v,v}$ and $\widetilde{T^v}$ (as mentioned in Proposition~\ref{prop.4.1} above) is a bounded RSK correspondence.
\subsection{Proof of the application}\label{ss.reduction}
Consider $v$ as an element of $I(d,2d)$. A \textbf{standard monomial} in $I(d,2d)$ is a totally ordered sequence $\theta_1\geq\ldots\geq\theta_t$ of elements of $I(d,2d)$. Such a monomial is called $v$-\textbf{compatible} if each $\theta_j$ is comparable to $v$ but no $\theta_j$ equals $v$; it is \textbf{anti-dominated} by $v$ if $\theta_t\geq v$. Let $\widetilde{SM^{v,v}}$ denote the set of all $v$-compatible standard monomials in $I(d,2d)$ anti-dominated by $v$.\par  
There is a natural injection $f:SM^{v,v}\rightarrow\widetilde{SM^{v,v}}$ given by 
$$f(\mathfrak{w}_1 \geq \ldots \geq \mathfrak{w}_t):=top(\mathfrak{w}_1)\geq bot(\mathfrak{w}_1)\geq\ldots\geq top(\mathfrak{w}_t)\geq bot(\mathfrak{w}_t).$$
Composing this map $f$ with the bijection $\tilde{\phi}$ (it is the inverse map of $\tilde{\pi}$ given in \S 4, \cite{kr}) from $\widetilde{SM^{v,v}}\rightarrow\widetilde{T^v}$, we get an injection of $SM^{v,v}$ into $\widetilde{T^v}$. It then follows from Lemma 4.5 of \cite{gr} that under this composition, the image of $SM^{v,v}$ in $\widetilde{T^v}$ is the set $\mathfrak{E}$ of all \textit{special monomials}. On the other hand, there is a bijective map (call it $g$) from the set $\mathfrak{E}$ of all special monomials to $\widetilde{T^v}$ as given in \S 4.1 of \cite{gr}: Given any $\mathfrak{S}$ in $\mathfrak{E}$, to get $g(\mathfrak{S})$, replace those $(r,c)$ of $\mathfrak{S}$ with $r>c^*$ by $(c^*,r^*)$ and then take the (positive) square root. The composition $\eta:=g\circ\tilde{\phi}\circ f$ is the required bijection from $SM^{v,v}$ to $\widetilde{T^v}$.\par 
Therefore, to prove that the composition map $\eta$ is a bounded RSK correspondence, it suffices to show that the map $\tilde{\phi}$ (of \S 4, \cite{kr}) is a bounded RSK correspondence. But the maps $\tilde{\pi}$ and $\tilde{\phi}$ (of \cite{kr}) are inverses of each other. Hence, it now suffices to show that the map $\tilde{\pi}$ of \cite{kr} is equal to the map $BRSK$ of \cite{kv}. This fact has been proved in Corollary \ref{c.main} above. 

\end{document}